\pdfoutput=1

\documentclass[english]{jnsao}
\usepackage[utf8]{inputenc}
\usepackage{lineno}

\usepackage[nameinlink,capitalize]{cleveref}

\usepackage{enumitem}
\usepackage{booktabs}
\usepackage{multirow}
\usepackage{graphicx}
\usepackage{mathtools}
\usepackage{algorithm,algpseudocode}
\usepackage{xcolor}
\usepackage{url}
\usepackage{algorithm}
\usepackage{algorithmicx}
\usepackage{footnote}
\usepackage{amsthm}
\usepackage{refcount}

\makesavenoteenv{tabular}
\makesavenoteenv{table}

\DeclareMathOperator*{\TV}{TV}

\DeclareMathOperator*{\argmax}{arg\,max}
\DeclareMathOperator*{\argmin}{arg\,min}

\newcommand*\dd{\mathop{}\!\mathrm{d}}
\newcommand{\R}{\mathbb{R}}
\newcommand{\N}{\mathbb{N}}

\newcommand{\B}{\mathcal{B}}
\newcommand{\symdif}{\vartriangle}

\theoremstyle{definition}
\newtheorem{assumption}{Assumption}[section]

\title{On convergence of binary trust-region steepest descent}

\author{Paul Manns\thanks{Faculty of Mathematics,
		TU Dortmund University (\url{paul.manns@tu-dortmund.de}).}
	\and 
	Mirko Hahn\thanks{Faculty of Mathematics,
		Otto von Guericke University Magdeburg.}
	\and
	Christian Kirches\thanks{Carl-Friedrich-Gau\ss -Faculty,
		Technical University of Braunschweig.}
	\and
	Sven Leyffer\thanks{Mathematics and Computer Science Division,
		Argonne National Laboratory.}
	\and 
	Sebastian Sager\footnotemark[2]}
\shortauthor{Manns, Hahn, Kirches, Leyffer, Sager}

\acknowledgements{This work was supported by the U.S.~Department of Energy, Office of Science, Office of Advanced Scientific Computing Research, Scientific Discovery 
through the Advanced Computing (SciDAC) Program through the FASTMath Institute under Contract No. DE-AC02-06CH11357 and from the German Research Foundation under GRK 2297 MathCoRe (project No. 314838170) and SPP~1962 (projects No. SA 2016/1-2, KI 1839/1-2) and from the German Federal Ministry of Education and Research within the program ``Mathematics for Innovations'' under the project ``Power to Chemicals.''
Argonne National Laboratory Preprint ANL/MCS-P9652-0222.}

\manuscriptsubmitted{2022-10-18}
\manuscriptaccepted{2023-07-19}
\manuscriptvolume{4}
\manuscriptnumber{10164}
\manuscriptyear{2023}
\manuscriptdoi{10.46298/jnsao-2023-10164}
\manuscriptcopyright{{\copyright} the authors}
\manuscriptlicense{CC-BY-NC-ND 4.0}

\begin{document}
\maketitle
\begin{abstract}
Binary trust-region steepest descent (BTR) and
combinatorial integral approximation (CIA) 
are two recently investigated approaches for the solution of
optimization problems with distributed binary-/discrete-valued
variables (control functions). We show improved convergence results for
BTR by imposing a compactness assumption that is similar to the
convergence theory of CIA. As a corollary we conclude that BTR also 
constitutes a descent algorithm on the continuous relaxation
and its iterates converge weakly-$^*$ to stationary points of
the latter. We provide computational results that validate our
findings. In addition, we observe a regularizing effect of
BTR, which we explore by means of a hybridization of
CIA and BTR.\\[2ex]
\\[2ex]
\end{abstract}

\section{Introduction}
\label{sec:introduction}

For bounded domains $\Omega \subseteq \R^d$ we are interested in optimization problems of the form
\begin{gather}\label{eq:p}
\inf_x J(x)\quad\text{s.t.}\quad x(s) \in \{0,1\} \text{ for almost all (a.a.) } s \in \Omega \text{ and } x \in L^2(\Omega),
\tag{P}
\end{gather}
$J$ is a map from $L^2(\Omega)$ to $\R$. For this problem class we study solutions of corresponding continuous relaxations of the form
\begin{gather}\label{eq:r}
\min_x J(x)\quad\text{ s.t. }\quad x(s) \in [0,1] \text{ for a.a.\ } s \in \Omega \text{ and } x \in L^2(\Omega)
\tag{R}
\end{gather}
and their relation to problem \eqref{eq:p}. The $\inf$ in the
formulation of \eqref{eq:p} and the $\min$ in the formulation of \eqref{eq:r} are 
deliberately
chosen to highlight that problem \eqref{eq:r} (in contrast to \eqref{eq:p})
admits a minimizer under mild assumptions
\cite{kirches2020compactness,manns2020multidimensional}.
We restrict ourselves to this setting in the interest of a concise 
presentation. One can, however, extend our analysis in different directions, for example, to the case where an $L^1$-regularization
term is added to the objective function, $J(x)$, by following the arguments in \cite{manns2021relaxed}. 

A rich class of instances of \eqref{eq:p} are mixed-integer PDE-constrained 
optimization problems, where $J = j \circ S$, where $j$ is the objective of the
optimization and $S$ the control-to-state operator of an underlying
partial differential equation (PDE). Such problems arise in many different areas such as topology optimization \cite{leyffer2021convergence,haslinger2015topology},
optimum experimental design \cite{yu2021multidimensional},
and gas network optimization \cite{hante2020mixed,hante2017challenges}.

We build on recent work on two algorithmic solution approaches of
\eqref{eq:p}.
Specifically we use insights of the available analysis
of combinatorial integral approximation (CIA) to improve the known
convergence results for \emph{binary trust-region steepest 
	descent} (BTR). We note that there
are more methods for (approximately) solving problems of
the form \eqref{eq:p}. For example, an established method for topology
optimization is the SIMP method \cite{bendsoe2004extensions},
which employs a non-convex penalization and thus regularization of
the controls (designs), and postprocessing of the solution of 
\eqref{eq:r} by means of lumping or filtering techniques.

\paragraph{Combinatorial integral approximation \cite{hante2013relaxation,kirches2020compactness,sager2012integer,sager2011combinatorial}}
The idea that underlies CIA is to split
the solution process of \eqref{eq:p} into solving the continuous
relaxation \eqref{eq:r} and then computing a $\{0,1\}$-valued approximation of the relaxed
solution. The approximation process can be analyzed in the weak-$^*$ topology of $L^\infty(\Omega)$
\cite{manns2020improved}. Relying on compactness properties of an
underlying control-to-state operator (e.g., $J = j \circ S$, where $j \in C(L^2(\Omega),\R)$ and $S : L^2(\Omega) \to L^2(\Omega)$ is a compact operator), a tight approximation of the optimal objective value by
the resulting approximants can be proved \cite{kirches2020compactness}. If stationary points
are computed in the first step, the approximation properties generalize accordingly.
We note that CIA can handle problem formulations, where the constraint $x(s) \in \{0,1\}$ is generalized
to $x(s) \in V$, where $V \subset \R$ is a finite set, by means of so-called special ordered set of type 1
(SOS1) reformulation of $V$ so that its elements become the vertices of a unit simplex in $\R^{|V|}$.

\paragraph{Binary trust-region steepest descent \cite{hahn2022binary,sharma2020inversion,vogt2020mixed}}
The BTR method solves trust-region subproblems in which the level sets of the control $x$,
corresponding to the values $0$ and $1$, are manipulated to greedily improve the linearized
objective. A trust-region constraint limits the volume of the level sets, which is the $L^1$-norm of the
control function, and can change from one accepted iterate to the next.
The analysis in \cite{hahn2022binary} shows
that BTR iterates eventually satisfy a condition called \emph{$\varepsilon$-stationarity}
under a regularity assumption used to obtain sufficient decrease of the
aggregated volume of level set manipulations.
Regarding the problem formulation \eqref{eq:p}, we note that generalizations of the BTR
method to the constraint $x(s) \in \{0,1\}^k$ for $k \in \N$ are conceivable but have not
been considered in the literature so far.

The concept \emph{$\varepsilon$-stationarity} as introduced in 
\cite{hahn2022binary} measures the projected
gradient of the objective and also provides, as we will show, a criticality measure for
first-order necessary
optimality conditions for \eqref{eq:r}.
Moreover, similar to CIA, the structural assumptions on the quantities that
appear in \eqref{eq:p} are in general not able to prevent 
a fine microstructure
from developing over the iterations. In fact, the weak-$^*$ closure of the feasible
set of \eqref{eq:p} in $L^\infty(\Omega)$ is the feasible set of its relaxation
\eqref{eq:r} \cite{lindenstrauss1966short,lyapunov1940completely}. However, it is not known whether the
iterates generated by BTR converge to
a limit point that satisfies a first-order
optimality condition of \eqref{eq:r}, in other words, 
if the termination tolerance of BTR is driven
to zero. This lack of a convergence result
is in contrast to CIA, which has stationary 
limits under a suitable compactness
assumption.

\paragraph{Standing assumptions}
We provide our comparison of CIA and BTR and
prove the new convergence result under the
following set of assumptions, which will
be discussed in detail in the remainder.
\begin{assumption}\label{ass:tr_convergence}~\\[-1em]
	\begin{enumerate}[label=(\alph*)]
		\item\label{itm:Jboundedbelow} Let $J : L^2(\Omega) \to \R$ be bounded from below.
		\item\label{itm:Jdifferentiable} Let $J : L^1(\Omega) \to \R$
		be Fr\'{e}chet differentiable.
		\item\label{itm:Jprime_Lipschitz}
		Let $\nabla J : L^1 (\Omega) \to L^2(\Omega)$ be Lipschitz continuous.
		\item\label{itm:Jprime_cc} Let $\nabla J : L^2(\Omega) \to L^2(\Omega)$ be
		weak-norm continuous (completely continuous).
	\end{enumerate}
\end{assumption}
\begin{remark}
	We note that the domain of $J$ can be restricted to $L^\infty(\Omega)$
	or the feasible set of \eqref{eq:r} for all considerations in this work. However, we
	require the continuity (differentiability) properties with respect to the $L^1(\Omega)$-norm on the domain
	space in \ref{itm:Jdifferentiable} and \ref{itm:Jprime_cc} in the remainder.
	We note that assuming continuity with respect to the codomain
	in $L^2(\Omega)$ is well defined for our purpose because all feasible points and
	iterates are also $L^\infty$-functions.
\end{remark}

\paragraph{Contributions}
We close the aforementioned theoretical gap between BTR and CIA. In particular, we
use the compactness 
\cref{ass:tr_convergence}, \ref{itm:Jprime_cc}
on the derivative of $J$, and we show that the BTR iterates produced by Algorithm 2
in \cite{hahn2022binary} converge weakly-$^*$ in  $L^\infty(\Omega)$
to a point that is feasible and satisfies a first-order optimality condition for the continuous
relaxation \eqref{eq:r}.
We perform several computational experiments on an example problem
that is governed by an elliptic PDE to validate our theoretical findings:
specifically,  BTR validates the near-optimality of the solution produced
by CIA. 

We have observed that BTR tends to produce controls whose
level sets have shorter interface lengths between them
in practice when started from zero or a thresholded 
control. At the same time it is able to produce 
objective values of similar quality as CIA.
While we cannot prove guarantees on this behavior, it 
motivates us to explore a hybrid method,
where we apply CIA but use a coarser control mesh in
order to compute the binary-valued approximation of the 
continuous relaxation. Then we start BTR from there, which 
allows us to combine the bounds and efficient 
running time behavior obtained with the CIA method while 
capitalizing on the regularization effect of BTR.

\paragraph{Structure of the paper}

In \S\ref{sec:cia} we formally introduce
the CIA method and show that its underlying 
approximation results hold under \cref{ass:tr_convergence}.
In \S\ref{sec:btsd-algo} we formally introduce and describe the BTR algorithm. In \S\ref{sec:convergence-BTSD-1Opoints} we relate it to \cite{hahn2022binary}
and state our main convergence result. The proof is presented in \S\ref{sec:proofs}.
We provide a computational validation of our findings, demonstrate the
aforementioned regularization effect, and investigate the observed regularization effects with a hybrid method
in \S\ref{sec:compu-valid}.
We provide auxiliary results in \S\ref{sec:appendix} and provide a brief discussion
of \cref{ass:tr_convergence} with respect to
the assumptions imposed in the earlier work 
\cite{hahn2022binary} in \S\ref{sec:relationship}.

\paragraph{Notation}

Let $d \in \N$ denote a dimension. For a measurable set $A \subset \Omega$, $\lambda(A)$ denotes the
Borel--Lebesgue measure of $A$ in $\R^d$. The function $\chi_A$ denotes the
$\{0,1\}$-valued characteristic function of the set $A$. Let $\B$ denote the
Borel $\sigma$-algebra on $\Omega$. For a set $A \subset \Omega$, the set $A^c$
denotes its complement in $\Omega$. For sets $A$, $B \subset \Omega$, the expression
$A \symdif B$ denotes the symmetric difference between $A$ and $B$,
that is, $A \symdif B := (A \cup B) \backslash (A \cap B)$.
The inner product of the Hilbert space $L^2(\Omega)$ is denoted by
$(\cdot,\cdot)_{L^2}$. For a space $X$ and its topological dual $X^*$ we denote the
pairing that puts $X$ and $X^*$ in duality by $\langle\cdot,\cdot\rangle_{X^*,X}$. We denote weak convergence with the
arrow $\rightharpoonup$ and weak-$^*$ convergence with the
arrow $\rightharpoonup^*$.

\section{Combinatorial Integral Approximation}\label{sec:cia}

CIA decomposes the solution process of \eqref{eq:p} into
two steps. First, the continuous relaxation \eqref{eq:r}
is solved (appproximately) and then the result is used 
to compute a sequence of $\{0,1\}$-valued functions
that are feasible for \eqref{eq:p} and converge to
the computed solution (or stationary point) of 
\eqref{eq:r} in the weak-$^*$ topology of 
$L^\infty(\Omega)$. The CIA algorithm is given in 
\cref{alg:cia}. Its key ingredients and asymptotics
are described below.
\begin{algorithm}[H]
	\caption{CIA algorithm to optimize \eqref{eq:p} and \eqref{eq:r}.}\label{alg:cia}
	\begin{flushleft}
		\textbf{Input:}
		$J : L^2(\Omega) \to \R$,
		$\nabla J: L^2(\Omega) \to L^2(\Omega)$.
		
		\textbf{Input:}
		Order-conserving domain dissection
		$\left(\mathcal{S}^n\right)_n \subset 2^{\mathcal{B}(\Omega)}$ (see \cref{dfn:order-conserving-domain-dissection}).
	\end{flushleft}
	\begin{algorithmic}[1]
		\State\label{ln:rel} $y \gets $ (Approximately) compute
		a stationary point of \eqref{eq:r}.
		\For{$n = 0,1,2,\ldots$}
		\State\label{ln:round} 
		$x^n \gets $
		\texttt{Round}$(y, \mathcal{S}^n)$.
		\EndFor
	\end{algorithmic}
\end{algorithm}

\paragraph{Inputs of \cref{alg:cia}}
\Cref{alg:cia} generally requires the objective 
function  $J$ as well as its gradient $\nabla J$ as 
inputs in order to solve \eqref{eq:r} in Line 
\ref{ln:rel}. We note that depending on the properties 
of \eqref{eq:r} and the chosen algorithm, this may be 
relaxed or strengthened and subgradients of $J$ may 
suffice (e.g., for projected subgradient methods)
or Hessian evaluations  (e.g., for a semi-smooth 
Newton's method) may be desirable.

The second step of CIA, that is the for-loop in
\cref{alg:cia}, requires a sequence of grids
$\left(\mathcal{S}^n\right)_n$ that decompose
the domain $\Omega$. These grids need to abide
a certain regularity that is defined in
\cref{dfn:order-conserving-domain-dissection}
in order to obtain the
aforementioned weak-$^*$ convergence, which
is explained in more detail below.

\paragraph{The subroutine \texttt{Round}
	and order-conserving domain dissections}
The subroutine \texttt{Round} in
Line \ref{ln:round} takes an
$L^\infty(\Omega)$-function $y$ that is $[0,1]$-valued
and a partition $\mathcal{S}^n$ of $\Omega$ as inputs and
computes a $\{0,1\}$-valued function $x^n$ from them.

In the literature on CIA, \eqref{eq:p} is usually 
formulated with finite sets $V \subset \R^m$, $m \in \N$, instead of $\{0,1\}$
as the co-domain of the optimization 
variables in \eqref{eq:p} and the output of Line \ref{ln:rel} is a convex 
coefficient function $\alpha : \Omega \to [0,1]^M$ that satisfies
$\sum_{i=1}^m \alpha_i(s) = 1$ a.e.
The rounding algorithm then transforms 
$\alpha$ to a function $\omega^n : \Omega \to \{0,1\}^M$
such that exactly one entry of $\omega^n(s)$ is one and all
others are zero a.e. To relate this to our setting, we 
can simply choose $m = 2$, $\alpha = (y, 1-y)^T$ and 
recover $x^n$ as $x^n = \omega_1^n$.

Under the assumption that the sequence of partitions
is an order-conserving domain dissection and with a suitable
implementation of \texttt{Round} that satisfies the prerequisites of
\cite[Proposition 3.5]{kirches2020compactness},
one obtains
\[ \alpha \rightharpoonup^* \omega \text{ in } L^\infty(\Omega), \]
which directly implies
\begin{gather}\label{eq:round_weakstar}
x^n \rightharpoonup^* y \text{ in } L^\infty(\Omega).
\end{gather}
Admissible choices for the subroutine \texttt{Round} are, for example,
sum-up rounding (SUR) \cite{sager2005numerical,manns2020multidimensional},
next-forced rounding \cite{jung2014relaxations},
and the combinatorial optimization-based algorithms in
\cite{bestehorn2019switching,jung2015lagrangian,zeile2021mixed}.
The key property of order-conserving domain dissections is that
during the refinement of the grid from one iteration to the next,
a spatial coherence property and a regular shrinkage property
that allow to leverage Lebesgue's differentiation theorem,
see the analysis in \cite{manns2020multidimensional}. A formal
definition is given in \cref{dfn:order-conserving-domain-dissection}
in \S\ref{sec:sur-appendix}. The choices for the subroutine \texttt{Round} that are used in our computational
experiments are described in \S\ref{sec:sur-appendix}.

\paragraph{Asymptotics
	of \cref{alg:binary_tr_steepest_descent}
	under \cref{ass:tr_convergence}}

Which set of assumptions is necessary so that solution algorithms
for \eqref{eq:r} produce sequences with (weak) cluster points that are stationary 
for \eqref{eq:r} in Line~\ref{ln:rel} depends on the properties of \eqref{eq:r} and the desired algorithm. 
If \eqref{eq:r} is convex, few assumptions may suffice and
Assumptions \ref{ass:tr_convergence}, \ref{itm:Jdifferentiable}, \ref{itm:Jprime_Lipschitz}, 
\ref{itm:Jprime_cc} may be relaxed to a boundedness assumption of the
$\varepsilon$-subdifferential on bounded sets and the projected (sub)gradient
method will work \cite{alber1998projected}.
If $J$ is not convex, a projected gradient method with standard
line search techniques, for example, Armijo linesearch, yields 
convergence to stationary points under
\cref{ass:tr_convergence} \cite{dunn1980convergence}. In this case,
a metricization of the domain space in \ref{itm:Jdifferentiable}
and \ref{itm:Jprime_Lipschitz} with the $L^2$-norm instead
of $L^1$-norm and the gradient needs not to be Lipschitz
but only uniformly continuous. More regularity allows to employ 
second-order methods like semi-smooth Newton to solve for
the first-order optimality conditions of \eqref{eq:r}.

While the assumption on the grids and the choice
of the \texttt{Round} subroutine imply \eqref{eq:round_weakstar},
the desired convergence of the objectives
\[ J(x^n) \to J(y) \]
in CIA requires that $J : L^2(\Omega) \to \R$ is weakly
continuous, which can often be asserted by regularity
(compactness) properties of an underlying differential
equation in the context of optimal control
\cite{kirches2020compactness,manns2020multidimensional}.
This is implied by Assumptions \ref{ass:tr_convergence}, 
\ref{itm:Jdifferentiable}, 
and \ref{itm:Jprime_cc} which is shown below.

\begin{proposition}\label{prp:J_weakcont}
	Let Assumptions \ref{ass:tr_convergence}, \ref{itm:Jdifferentiable}, 
	and \ref{itm:Jprime_cc} hold. Then $J : L^2(\Omega) \to \R$
	is weakly continuous.
\end{proposition}
\begin{proof}
	Let $x^n$, $x \in L^2(\Omega)$ be such that $x^n \rightharpoonup x$, meaning $x^n$
	converging weakly to $x$, in $L^2(\Omega)$.
	We need to show $J(x^n) \to J(x)$. \Cref{ass:tr_convergence}, \ref{itm:Jdifferentiable}
	and the mean value theorem imply that $J(x^n) - J(x) = (\nabla J(\xi^n), x^n - x)_{L^2}$
	for some $\xi^n \in L^2(\Omega)$ in the line segment between $x^n$ and $x$ for all $n \in \N$.
	Because $(\xi^n)_n$ is bounded, there exists a weakly convergent subsequence $\xi^{n_k} \rightharpoonup \xi$
	for some $\xi \in L^2(\Omega)$. \Cref{ass:tr_convergence}, \ref{itm:Jprime_cc}
	implies $\nabla J(\xi^{n_k}) \to \nabla J(\xi)$ in $L^2(\Omega)$, and consequently
	$J(x^{n_k}) - J(x) = (\nabla J(\xi^{n_k}), x^{n_k} - x)_{L^2} \to 0$.
	Passing to subsubsequences proves the claim.
\end{proof}

\begin{remark}
	We note that assumption of an order-conserving domain
	dissection and the weak continuity of $J$ are sufficient
	to obtain the desired weak-$^*$ convergence of the $x^n$
	to stationary points and the corresponding convergence
	objective values as well if one does not compute $y$ first and
	then executes \texttt{Round} but instead executes \texttt{Round}
	on the iterates of produced by an optimization algorithm
	for \eqref{eq:r}, see Theorem 4.7 in \cite{manns2020multidimensional}.
\end{remark}

\section{Binary Trust-Region Steepest Descent}
\label{sec:btsd-algo}

The BTR algorithm
operates on characteristic functions induced by measurable sets.
We introduce the inputs and the trust-region subproblem. Then we describe the iterations of 
\cref{alg:binary_tr_steepest_descent} step by step.
We relate the quantities in our variant of the algorithm to the one introduced as Algorithm 2
in \cite{hahn2022binary}, which purely takes the point of view of measurable sets.

The BTR algorithm is given as \cref{alg:binary_tr_steepest_descent} and 
is a special case of \cite[Algorithm\,2]{hahn2022binary}. In particular,
it corresponds to \cite[Algorithm 2]{hahn2022binary}
with the choices $\mathcal{J}(A) \coloneqq J(\chi_A)$ for $A \in \B$.
We also choose $\varepsilon = 0$ because we aim to
study the asymptotics of the algorithm when it is not stopped early.

\begin{algorithm}[H]
	\caption{BTR algorithm to optimize \eqref{eq:p} and
		\eqref{eq:r}.\label{alg:binary_tr_steepest_descent}}
	\begin{flushleft}
		\textbf{Input:}
		$J : L^2(\Omega) \to \R$,
		$\nabla J: L^2(\Omega) \to L^2(\Omega)$,
		$\Delta_{\max} \in (0,\lambda(\Omega))$,
		$0 < \sigma_1 < \sigma_2 \le 1$,
		$\omega \in (0,1)$.
		
		\textbf{Input:} $U^0 \in \B$, $\Delta^0 \in (0, \Delta_{\max})$
	\end{flushleft}
	\begin{algorithmic}[1]
		\For{$n = 0,1,2,\ldots$}
		\State\label{ln:gUn_approx}$g^n \gets \nabla J(\chi_{U^n})(\chi_{(U^n)^c} - \chi_{U^n})$
		\State\label{ln:find_step}$D^n \gets$ \texttt{FindStep} \big($g^n, \Delta^{n},
		\min\{\omega, 0.5 \|\min\{g^n, 0\}\|_{L^1}/\lambda(\Omega), \Delta^{n}\}$\big)
		\quad// e.g.\ \cite[Proc.\,1]{hahn2022binary}
		\If{$J(\chi_{U^n \symdif D^n}) - J(\chi_{U^n})
			\le \sigma_1 (\nabla J(\chi_{U^n}), \chi_{D^n\setminus U^n} - \chi_{U^n\cap D^n}
			)_{L^2}$}\label{ln:sufficient_decrease_condition}
		\State $U^{n+1} \gets U^{n} \symdif D^n$
		\If{$J(\chi_{U^n \symdif D^n}) - J(\chi_{U^n}) \le \sigma_2  (\nabla J(\chi_{U^n}), \chi_{D^n\setminus U^n} - \chi_{U^n\cap D^n})_{L^2}$}
		\State $\Delta^{n+1} \gets \min\{2\Delta^n, \Delta_{\max}\}$
		\Else
		\State $\Delta^{n+1} \gets \Delta^n$
		\EndIf
		\Else
		\State $(U^{n+1}, \Delta^{n+1}) \gets (U^{n}, 0.5 \Delta^{n})$
		\EndIf
		\EndFor
	\end{algorithmic}
\end{algorithm}

\paragraph{Inputs of \cref{alg:binary_tr_steepest_descent}}
The algorithm requires the objective function $J$ as well as its
gradient $\nabla J$ as inputs. 
Using the latter requires assuming differentiability of $J$ with respect to the $L^2$-norm.
For the acceptance criterion of the computed descent step and
the update of the trust-region radius, the algorithm requires a maximal trust-region radius 
$\Delta_{\max}$ and control parameters $\sigma_1$ and $\sigma_2$ as inputs.
To compute a descent step, the algorithm uses the subroutine \texttt{FindStep} (see below). The subroutine
also requires the parameter $\omega$, which ensures that the volume of
the returned set is always bounded from below by a fraction of $\Delta^n$
that is smaller than $1$.

\paragraph{Trust-region subproblem and subroutine
	\texttt{FindStep}}
The subroutine \texttt{FindStep} in Line \ref{ln:find_step} of
\cref{alg:binary_tr_steepest_descent} approximately solves the subproblem
\begin{gather}\label{eq:g_subproblem}
\min_{D} \int_D g(s)\dd s
\quad\text{ s.t.\ }\quad \left\{
\begin{aligned}D \subset g^{-1}((-\infty,0]),&\\
\lambda(D) \le \Delta.
\end{aligned}
\right.
\end{gather}
In \cref{alg:binary_tr_steepest_descent} 
\texttt{FindStep} is called with
$g = \nabla J(\chi_{U^n})(\chi_{(U^n)^c} - \chi_{U^n})$
and $\Delta = \Delta^n$. Changing from set optimization
to function optimization, the minimization
problem \eqref{eq:g_subproblem} is equivalent to
the minimization problem
\begin{gather}\label{eq:tr}
\min_{d}\ (\nabla J(\chi_{U}),d)_{L^2} \quad\text{s.t.}\quad 
\left\{
\begin{aligned}
\chi_{U}(s) + d(s) \in \{0,1\} &\text{ for a.a.\ } s \in \Omega, \\
\|d\|_{L^1} \le \Delta
\end{aligned}
\right.
\tag{TR(\protect\ensuremath{\Delta})}
\end{gather}
if $g = \nabla J(\chi_{U})(\chi_{U^c} - \chi_{U})$. We provide a proof of the equivalence in \cref{prp:tr_equivalence}. We note that
using $\min$ in the definitions of 
\eqref{eq:g_subproblem} and \eqref{eq:tr}
is justified because \eqref{eq:tr}, and thus also \eqref{eq:g_subproblem}, indeed admits a 
minimizer, which is shown in 
\cref{prp:min_tr_subproblem}.

The analysis in \cite{hahn2022binary}
employs that, in every iteration, 
\texttt{FindStep} produces a set $D$,
or, equivalently, 
a corresponding function $d = \chi_{U\symdif D} - \chi_{U}$ (see also \cref{prp:tr_equivalence}),
such that
\begin{align*}
(\nabla J (\chi_{U}), d)_{L^2}
&\le \min \left\{ (\nabla J (\chi_{U}), s)_{L^2}
\,\middle|\, s \text{ feasible for \eqref{eq:tr}} \right\}
+ \delta \Delta \ \text{ and}\\ 
\|d\|_{L^1} &\le \Delta
\end{align*}
hold, where
$\delta = \min\left\{\omega, 0.5\|\min\{g, 0\}\|_{L^1}/\lambda(\Omega), \Delta\right\}$
is the third input of \texttt{FindStep} in 
\cref{alg:binary_tr_steepest_descent}.
A bisection algorithm in function space that realizes this property
is also provided in \cite{hahn2022binary}.
Other algorithmic approaches to (approximately)
solving \eqref{eq:tr} and thus implementing \texttt{FindStep} are
possible as well. In \S\ref{sec:compu-valid} we describe
and take advantage of a variant that exploits uniform  meshes for the 
control discretization in our computational experiments.

\paragraph{Description of the steps of \cref{alg:binary_tr_steepest_descent}}
The for-loop
starting in Line 1 of \cref{alg:binary_tr_steepest_descent} computes candidates for improvements
of the objective function that are either accepted or rejected and then updates the trust-region 
radius accordingly.

Line \ref{ln:gUn_approx} computes the function $g^n$ such that it is equal to
$-\nabla J(\chi_{U^n})$ on $U^n$, where $\chi_{U^n}$ may be decreased, and such that it is equal to
$\nabla J(\chi_{U^n})$ on $(U^n)^c$, where $\chi_{U^n}$ may be increased. Thus
$\{s \in \Omega\,|\,g^n(s) \le 0\}$ is the set on which $U^n$ can be changed to obtain
a first-order decrease of $J$. Because of the use of the $L^1$-norm, the trust-region radius 
$\Delta^n$ limits the volume of the set that can be changed in the current iteration, see
\cref{prp:tr_equivalence}.

The subroutine \texttt{FindStep} in Line \ref{ln:find_step} computes a set
$D^n \subset (g^n)^{-1}((-\infty,0])$ that 
approximates the solution of~\eqref{eq:g_subproblem}.

The candidate for improving $J$ is then the modification of the characteristic function of
the set $U^n$, where the values on $D^n$ are flipped, or formally $\chi_{U^n \symdif D^n}$. 
From Line \ref{ln:sufficient_decrease_condition} onward, the for-loop resembles common 
trust-region methods. Line \ref{ln:sufficient_decrease_condition} determines whether the 
reduction achieved by $\chi_{U^n \symdif D^n}$ is at least a fraction of the reduction predicted by the linear model, in which case the step is accepted.
A second (larger) ratio is used to determine whether the trust-region radius should be increased
(doubled) or left unchanged after acceptance. The trust-region radius is reduced (halved) after rejection of a candidate step.

\section{Convergence of BTR to First-Order Optimal Points}
\label{sec:convergence-BTSD-1Opoints}

We introduce our main result. \cref{alg:binary_tr_steepest_descent}
operates on iterates that are feasible
for the integer problem \eqref{eq:p}. We prove
that our algorithm generates a sequence of
integer feasible points whose limits
are first-order optimal points of the relaxation \eqref{eq:r}. Thus,
we obtain a minimizing sequence for \eqref{eq:p}, which itself may not have a solution,
if \eqref{eq:r} is a convex problem.
We introduce and a criticality measure
for \eqref{eq:r} and relate it to \cref{alg:binary_tr_steepest_descent}
before introducing the main theorem.

\subsection{Criticality Measure for \eqref{eq:r}}
\label{sec:criticality_measure}

We define the criticality measure $C : L^2(\Omega) \to [0,\infty)$ for $x \in L^2(\Omega)$ as
\begin{gather}\label{eq:criticality}
C(x) \coloneqq
\max \left\{
\int_\Omega \nabla J(x)(x - f)\dd s\,\middle|\,
f \text{ feasible for } \eqref{eq:r}\right\}
=   \int_\Omega \nabla J(x)x
+ \int_\Omega \max\{-\nabla J(x), 0\}\dd s,
\end{gather}
where the identity follows from the structure of the
feasible set of \eqref{eq:r}.
$C$ coincides with the function $\Phi$ in 
\cite{dunn1980convergence} and is also known as
\emph{primal gap function} \cite{larsson1994class}.
Local minimizers are zeros of $C$, which is well known
and repeated here for convenience. In particular,
this leads to the usual definition of stationary points
below.
\begin{proposition}\label{prp:local_optimality_condition}
	If $x$ is a local minimizer of \eqref{eq:r},
	then $\nabla J(x)(x - f) \le 0$ a.e.\ holds for all $f$ that
	are feasible for \eqref{eq:r}. Moreover, $C(x) = 0$.
\end{proposition}
\begin{proof}
	The first claim follows from a Taylor expansion at $x$.
	$C(x) = 0$ because $x$ is feasible in the $\max$.
\end{proof}
\begin{definition}
	A function $x$ that is feasible for \eqref{eq:r}
	is called \emph{stationary} if $C(x) = 0$.
\end{definition}
The \texttt{FindStep} subroutine in 
\cref{alg:binary_tr_steepest_descent} operates
with the quantity $\|\min\{g,0\}\|_{L^1}$ with the
choice $g = \nabla J(\chi_U)(\chi_{U^c} - \chi_U)$
for a set $U \in \B$. We show below that
$C(\chi_U) = \|\min\{g,0\}\|_{L^1}$. In this case,
the first-order optimality condition 
$C(\chi_U) = 0$ from \cref{prp:local_optimality_condition}
corresponds the first-order optimality condition
of the set-based view point in
\cite{hahn2022binary}, see Lemma 5 and Corollary 1
therein.

\begin{proposition}\label{prp:criticality_char}
	Let $U \in \B$. Then for
	$g = \nabla J(\chi_U)(\chi_{U^c} - \chi_U)$ it holds that
	$C(\chi_U) = \|\min\{g,0\}\|_{L^1}$.
\end{proposition}
\begin{proof}
	For a.a.\ $s \in\Omega$, we obtain
	\begin{gather}\label{eq:criticality_characterization}
	\min\{g(s), 0\} = \left\{
	\begin{aligned}
	-\nabla J(\chi_{U})(s) &\ \text{ if }s \in U \text{ and } \nabla J(\chi_{U})(s) \ge 0,\\
	\nabla J(\chi_{U})(s) &\ \text{ if }s \in {U}^c \text{ and } \nabla J(\chi_{U})(s) \le 0,\\
	0 &\ \text{ else.}
	\end{aligned}
	\right.
	\end{gather}
	For the choice $x = \chi_U$ in
	\eqref{eq:criticality}, the right hand side implies
	\[ C(x) = \int_\Omega \left\{
	\begin{aligned} |\nabla J(x)(s)| &\ \text{ if }
	(s \in U \text{ and } \nabla J(x)(s) \ge 0)
	\text{ or }
	(s \in U^c \text{ and } \nabla J(x)(s) \le 0),\\
	0 &\ \text{ else}
	\end{aligned}\right\} \dd s. \]
	The claim follows by a pointwise a.e.\ comparison of the integrands.
\end{proof}
An alternative criticality measure for \eqref{eq:r} that uses the $L^1$-norm is the
function $\tilde{C}$ that is defined for
$x \in L^1(\Omega)$ as
\[ \tilde{C}(x) \coloneqq \|x - P_{[0,1]}(x - \nabla J(x))\|_{L^1}.\]
It is known that $\tilde{C}(x) = 0$ implies that $x$ 
satisfies a  first-order necessary optimality condition
for \eqref{eq:r}; see, for example,
\cite[Lemma\,1.12]{hinze2008optimization}.
Moreover, it can easily be verified that
$C(x) = 0$ if and only if $\tilde{C}(x) = 0$.
We do not use $\tilde{C}$ because it would complicate our analysis and lead
to a less concise presentation.

\subsection{Main Result}
\label{sec:main_result}
Having introduced the necessary notation, concepts, and assumptions, we  now  state
our main convergence results.

\begin{theorem}\label{thm:asymptotics_part_1}
	Let Assumptions \ref{ass:tr_convergence}, \ref{itm:Jboundedbelow} and \ref{itm:Jdifferentiable}
	hold. Let $(U^n)_n \subset \B$, $(D^n)_n \subset \B$, and
	$(\Delta^n)_n \subset (0,\Delta_{\max}]$
	be the sequences of sets and trust-region radii produced by \cref{alg:binary_tr_steepest_descent}.
	Then the sequence of objective values $(J(\chi_{U^n}))_n$ is monotonically nonincreasing.
	Moreover, one of the following mutually exclusive outcomes holds:
	\begin{enumerate}
		\item\label{itm:outcome_finite_stationary} There exists $n_0 \in \N$ such that $U^{n_0} = U^{n}$ a.e.\ holds for all $n \ge n_0$. Then $\chi_{U^{n_0}}$ is stationary for \eqref{eq:r}.
		\item\label{itm:outcome_accumulation_stationary} For all $n_0 \in\N$ there exists $n_1 > n_0$ such that $\lambda(U^{n_1} \symdif U^{n_0}) > 0$.
		The sequence $(\chi_{U^n})_n \subset L^\infty(\Omega)$ admits a weak-$^*$ accumulation point.
		Every weak-$^*$ accumulation point $f$ of $(\chi_{U^n})_n$ is feasible for \eqref{eq:r}.
		
		If additionally \cref{ass:tr_convergence}, \ref{itm:Jprime_Lipschitz} holds, 
		then 
		\[ \lim_{n\to\infty} C(\chi_{U^n}) = 0. \]
		If additionally \cref{ass:tr_convergence}, \ref{itm:Jprime_cc} holds  and if a subsequence $(\chi_{U^{n_k}})_k \subset (\chi_{U^n})_n$ satisfies
		$C(\chi_{U^{n_k}}) \to 0$, then 
		every weak-$^*$ accumulation point $f$ of $(\chi_{U^{n_k}})_k$ satisfies
		$C(f) = 0$, i.e., $f$ is stationary for \eqref{eq:r}.
	\end{enumerate}
\end{theorem}

\Cref{thm:asymptotics_part_1} is proven
in \S\ref{sec:proofs}.
We obtain the following corollary that
shows that \cref{alg:binary_tr_steepest_descent} produces a sequence
of binary iterates that converge weakly-$^*$ to stationary points of the
continuous relaxation \eqref{eq:r} of \eqref{eq:p}. Thus BTR
yields results comparable to those produced by CIA because
solution algorithms for \eqref{eq:r} cannot be expected to perform better than
producing a stationary point of \eqref{eq:r} in practice.

\begin{corollary}
	Let \cref{ass:tr_convergence} hold.
	Let $(U^n)_n \subset \B$ be the sequences of sets produced by
	Algorithm \ref{alg:binary_tr_steepest_descent}.
	Then all weak-$^*$ accumulation points of $(\chi_{U^n})_n$ are feasible and stationary
	for \eqref{eq:r}.
\end{corollary}
\begin{proof}
	The claim follows by combining the two assertions
	in Outcome \ref{itm:outcome_accumulation_stationary}
	of \cref{thm:asymptotics_part_1}.
\end{proof}

\section{Proof of the Main Theorems}
\label{sec:proofs}

In this section we prove \cref{thm:asymptotics_part_1}. 
We first prove preparatory results on the sufficient reduction condition
(\cref{alg:binary_tr_steepest_descent}, Line~\ref{ln:sufficient_decrease_condition})
for binary-valued control functions and trust-region steps in \S\ref{sec:sufficient_decrease}.
Then we employ these results to analyze the asymptotics of
\cref{alg:binary_tr_steepest_descent} in \S\ref{sec:proof_asymptotics},
finishing with the proof of \cref{thm:asymptotics_part_1}.

\subsection{Sufficient Decrease with a Characteristic Function}\label{sec:sufficient_decrease}

The first step of the proof is to show that if $\chi_{U^n}$ for
a given iterate $U^n \in \B$ of \cref{alg:binary_tr_steepest_descent}
is not stationary, then there exists a set $D^n$ such that $\chi_{U^n \symdif D^n}$ 
satisfies a sufficient decrease condition with respect to $\nabla J(\chi_{U^n})$ for 
sufficiently small trust-region radii.
We briefly recap and the well-known result on existence of a 
descent direction in  
\cref{lem:descent_binary_construction_var},
which we adapt for our case of characteristic
functions as \cref{alg:binary_tr_steepest_descent}
operates on. This in turn implies acceptance of a new
iterate after finitely many steps as is shown in \cref{lem:finite_step_acceptance}.
\begin{lemma}\label{lem:descent_binary_construction_var}
	Let \cref{ass:tr_convergence} \ref{itm:Jdifferentiable} 
	hold. Let $\chi_U$ for $U \in \B$ not be stationary for 
	\eqref{eq:r}. Then there exist $\varepsilon > 0$ and
	$\Delta_0 > 0$ such that for all
	$0 < \Delta \le \Delta_0$, there exists $d \in L^1(\Omega)$
	that is feasible for \eqref{eq:tr} and satisfies $(\nabla J(\chi_U), d) \le
	-\varepsilon\Delta$.
\end{lemma}
\begin{proof}
	Because $\chi_U$ is not stationary,
	we have $C(\chi_U) > 0$, which implies that there
	is $\varepsilon > 0$ and a set $D \subset g^{-1}((-\infty,0])$ with $\lambda(D) > 0$
	such that $\nabla J(\chi_U)(s)(\chi_{U^c}(s) - \chi_U(s)) < -\varepsilon$ for a.a.\ $s \in D$.
	Let $\Delta_0 \coloneqq \lambda(D)$ and using
	the regularity properties of the Lebesgue measure,
	there is a subset $D_\Delta \subset D$ with $\lambda(D_\Delta) = \Delta$ for all
	$0 < \Delta \le \Delta_0$ (note that it is possible
	to use the greedy construction from \cref{prp:min_tr_subproblem}
	here too), which implies
	\[ \int_{D_\Delta} \nabla J(\chi_U)(\chi_{U^c} - \chi_U)
	\dd s
	< -\Delta \varepsilon.\]
	Using the equivalence asserted in \cref{prp:tr_equivalence} and in particular setting
	$d \coloneqq \chi_{U \symdif D_\Delta} - \chi_U$
	yields the claim.
\end{proof}

We employ this result to prove that \cref{alg:binary_tr_steepest_descent}
accepts a step after finitely many iterations if the current iterate is not stationary
for \eqref{eq:r}.

\begin{lemma}\label{lem:finite_step_acceptance}
	Let Assumptions \ref{ass:tr_convergence}, \ref{itm:Jdifferentiable}
	and \ref{itm:Jprime_Lipschitz}
	be satisfied. Let $(U^n)_n \subset \B$, $(D^n)_n \subset \B$, and $(\Delta^n)_n \subset (0,\Delta_{\max}]$
	be the sequences of sets and trust-region radii produced by \cref{alg:binary_tr_steepest_descent}.
	Let $\chi_{U^n}$ not be stationary for \eqref{eq:r}. Then
	the output of \cref{alg:binary_tr_steepest_descent}, Line \ref{ln:find_step},
	is accepted after $k \in \N$ steps: specifically,
	$U^n = U^{n + j}$ for all $0 \le j < k$ and
	$J(\chi_{U^n \symdif D^{n + k}}) - J(\chi_{U^n}) \le \sigma_1 (\nabla J(\chi_{U^n}),\chi_{D^{n+k}\setminus U^n} - \chi_{U^n \cap D^n})_{L^2}$.
\end{lemma}
\begin{proof}
	For $m \in \N$, we define the optimal linear predicted reduction as
	\[ L^m \coloneqq -\inf_{D' \in \B}\left\{ (\nabla J(\chi_{U^m}),d')_{L^2} \,\Big|\,
	d' = \chi_{D'\setminus U^m} - \chi_{D' \cap U^m} \text{ and }
	\|d'\|_{L^1} \le \Delta^m
	\right\}. \]
	By design of \cref{alg:binary_tr_steepest_descent}, we have
	$\Delta^{n + j + 1} = 0.5 \Delta^{n + j}$ for all $j \ge 0$ until the step $D^{n+j}$ is accepted.
	
	We prove the claim by contradiction and assume that the step $D^{n+k}$ is not accepted for all $k \in \N$.
	Because $U^n$ is not stationary for \eqref{eq:r} and $\Delta^{n + k} \to 0$ for $k \to \infty$,
	\cref{lem:descent_binary_construction_var} implies that there exist $\varepsilon > 0$ and
	$k_0 \in \N$ such that for all $k \ge k_0$ the estimate $L^{n + k} \ge \Delta^{k + j}\varepsilon$ holds.
	
	We apply Taylor's theorem and obtain that
	\begin{align}
	\hspace{2em}&\hspace{-2em}J(\chi_{U^n}) - J(\chi_{U^n \symdif D^{n + k}})\nonumber \\
	&= -(\nabla J(\chi_{U^n}), \chi_{D^{n+k}\setminus U^n}
	- \chi_{U^n\cap D^{n+k}})_{L^2} + o\left(\lambda(D^{n+k})\right)\nonumber\\
	&\ge L^{n+k} - (\Delta^{n + k})^2 + o(\Delta^{n+k})\label{eq:2_latter_terms},
	\end{align}
	where the inequality follows from the construction of $D^{n+k}$ by means of
	the \texttt{FindStep} subroutine, specifically Lemma 9 in \cite{hahn2022binary} with a choice
	$\delta \le \Delta^{n+k}$. The inequality $\delta \le \Delta^{n+k}$ is satisfied
	in \cref{alg:binary_tr_steepest_descent},
	Line \ref{ln:find_step}, because the parameter $\delta$ of \texttt{FindStep}, the third argument
	of the subroutine, is given a value that is less than or equal to the trust-region radius in all
	iterations.
	
	Because $L^{n+k} \ge \varepsilon \Delta^{n+k}$ holds for all $k \ge k_0$ and the
	two latter terms in \eqref{eq:2_latter_terms} are $o(\Delta^{n+k})$
	there exists $k_1 \in \N$ such that for all $k \ge k_1$ it holds that
	\[ J(\chi_{U^n}) - J(\chi_{U^n \symdif D^{n + k}}) \ge \sigma_1 L^{n+k}. \]
	By definition of $L^{n+k}$ it follows that
	\[ J(\chi_{U^n}) - J(\chi_{U^n \symdif D^{n + k}})
	\ge \sigma_1 L^{n+k}
	\ge -\sigma_1 (\nabla J(\chi_{U^n}),\chi_{D^{n+k}\setminus U^n} - \chi_{U^n \cap D^{n+k}})_{L^2}, \]
	and thus the step $D^{n+k_1}$ is accepted in \cref{alg:binary_tr_steepest_descent}.
	This contradicts the assumption that the step $D^{n+k}$ is not accepted for all $k \in \N$.
\end{proof}

\subsection{Asymptotics of \cref{alg:binary_tr_steepest_descent}}\label{sec:proof_asymptotics}

Before finalizing the proof of \cref{thm:asymptotics_part_1},
we show three further preparatory lemmas. \Cref{lem:non_increasing_objectives} states that the
sequence of iterates produced by \cref{alg:binary_tr_steepest_descent} has a corresponding sequence
of monotononically nonincreasing objective values. \Cref{lem:criticality_coercion_implies_tr_contraction}
shows that if the criticality measure $C$ stays bounded away from zero over the
iterations of \cref{alg:binary_tr_steepest_descent}, then the trust-region radius contracts
to zero.

\begin{lemma}\label{lem:non_increasing_objectives}
	Let \cref{ass:tr_convergence}, \ref{itm:Jdifferentiable} hold.
	Let $(U^n)_n$, $(D^n)_n \subset \B$, and $(\Delta^n)_n \subset (0,\Delta_{\max}]$
	be the sequences of sets and trust-region radii produced by \cref{alg:binary_tr_steepest_descent}.
	Then the sequence of objective values $(J(\chi_{U^n}))_n$ is monotonically nonincreasing.
\end{lemma}
\begin{proof}
	By construction of $D^n$ with \texttt{FindStep}, Procedure 1 of \cite{hahn2022binary}, it holds that
	$D^n \subset \{s\in\Omega\,|\,g_{U^n}(s) < 0\}$. A step that is accepted in
	\cref{alg:binary_tr_steepest_descent} Line \ref{ln:sufficient_decrease_condition} satisfies
	$J(\chi_{U^{n} \symdif D^n}) < J(\chi_{U^n})$ because
	\[ (\nabla J(\chi_{U^n}), \chi_{D^n \setminus U^n} - \chi_{U^n\cap D^n})_{L^2}
	= \int_{D^n} g^n\,\dd s < 0, \]
	while $J(\chi_{U^n})$ remains unchanged for rejected steps.
	Thus, the sequence of objective values $(J(\chi_{U^n}))_n$ is monotonically nonincreasing.
\end{proof}

\begin{lemma}\label{lem:criticality_coercion_implies_tr_contraction}
	Let Assumptions \ref{ass:tr_convergence}, \ref{itm:Jboundedbelow} and \ref{itm:Jdifferentiable}
	hold.
	Let $(U^n)_n \subset \B$, $(D^n)_n \subset \B$, and $(\Delta^n)_n \subset (0,\Delta_{\max}]$
	be the sequences of sets and trust-region radii produced by \cref{alg:binary_tr_steepest_descent}.
	If there exists $\varepsilon > 0$ and $n_0 \in \N$
	such that $C(\chi_{U^n}) > \varepsilon$
	for all $n \ge n_0$, then $\Delta^n \to 0$.
\end{lemma}
\begin{proof}
	We use the notation $L^m$ for the optimal predicted reduction in iteration $m \in \N$ as in
	the proof of \cref{lem:finite_step_acceptance}.
	From \cref{prp:min_tr_subproblem} and the definition of
	$C$ it follows that
	$L^m \ge C(\chi_{U^m})\frac{\Delta^m}{\lambda(\Omega)}$ for all iterations
	$m \in \N$. This can be seen by using the greedily constructed set.
	From \cref{alg:binary_tr_steepest_descent} Line \ref{ln:find_step}, \cref{prp:criticality_char},
	and
	$\Delta_{\max} \le \lambda(\Omega)$ it follows that the third parameter $\delta$ of the subroutine
	\texttt{FindStep} satisfies $\delta \le C(\chi_{U^n})/(2\Delta^n)$.
	The analysis of \texttt{FindStep}, specifically \cite[Lemma 9]{hahn2022binary}, implies
	that
	\begin{align*}
	-(J(\chi_{U^n}),\chi_{D^n\setminus U^n} - \chi_{U^n\cap D^n})_{L^2}
	&\ge L^n - \min\{\omega, C(\chi_{U^n})/(2\lambda(\Omega)), \Delta^n\}\Delta^n\\
	&\ge C(\chi_{U^n})\frac{\Delta^n}{\lambda(\Omega)} - \frac{1}{2}C(\chi_{U^n})\frac{\Delta^n}{\lambda(\Omega)}
	\ge \frac{\varepsilon\Delta^n}{2\lambda(\Omega)}
	\end{align*}
	for all $n \in \N$ and thus all $n \ge n_0$.
	
	We close the proof with a contradictory argument and assume that $\Delta^n \not\to 0$. We deduce
	that there exists an infinite subsequence $(\Delta^{n_k})_k$ of $(\Delta^n)_n$ such that
	$\liminf_{k\to\infty} \Delta^{n_k} > \underline{\Delta}$ for some $\underline{\Delta} > 0$.
	Consequently, there exists an infinite subsequence $(n_\ell)_\ell$ of accepted iterates 
	with trust-region radii $\Delta^{n_\ell} \ge \underline{\Delta}$.
	
	Combining these insights on the accepted iterates with the lower bound on the linearly predicted
	reductions, we obtain that
	\[ J(\chi_{U^{n_\ell}}) - J(\chi_{U^{n_\ell}\symdif D^{n_\ell}})
	\ge \sigma\frac{\varepsilon\underline{\Delta}}{4\lambda(\Omega)}.
	\] 
	Because the sequence of objective values is monotonically 
	nonincreasing by virtue of
	\cref{lem:non_increasing_objectives}, this implies $J(\chi_{U^{n_\ell}}) \to -\infty$,
	which contradicts \cref{ass:tr_convergence}, \ref{itm:Jboundedbelow} and thus the 
	assumption that $\Delta^n \not\to 0$.
\end{proof}

We are now ready to finish the proof of the two
main results.

\begin{proof}[Proof of \cref{thm:asymptotics_part_1}]
	\Cref{lem:non_increasing_objectives} proves the claim that the sequence of objective values is monotonically nonincreasing.
	
	We first analyze Outcome \ref{itm:outcome_finite_stationary}. Because there exists $n_0 \in \N$ such that 
	$U^{n_0} = U^{n}$ holds a.e.\ for all $n \ge n_0$, the acceptance criterion in \cref{alg:binary_tr_steepest_descent},
	Line \ref{ln:sufficient_decrease_condition}, is violated for all $n \ge n_0$. Then the claim of Outcome 
	\ref{itm:outcome_finite_stationary} follows from \cref{lem:finite_step_acceptance}.
	
	If there is no $n_0 \in \N$ such that $U^{n_0} = U^{n}$ holds a.e.\ for all $n \ge n_0$, then Outcome 
	\ref{itm:outcome_finite_stationary} does not hold true, and for all $n_0 \in\N$ there exists $n_1 > n_0$
	such that $\lambda(U^{n_1} \symdif U^{n_0}) > 0$. It follows that Outcomes 
	\ref{itm:outcome_finite_stationary} and \ref{itm:outcome_accumulation_stationary} are mutually exclusive.
	Moreover, the sequence $(\chi_{U^n})_n \subset L^\infty(\Omega)$ is bounded and thus admits a
	weak-$^*$ cluster point. By virtue of, for example, \cite[Theorem\,3]{tartar1979compensated},
	every weak-$^*$ cluster point of $(\chi_{U^n})_n$ is feasible for \eqref{eq:r}.
	
	Next we assume that Assumptions 
	\ref{ass:tr_convergence}, \ref{itm:Jboundedbelow}, \ref{itm:Jdifferentiable},
	and \ref{itm:Jprime_Lipschitz}
	are satisfied and prove the claim
	$\lim_{n\to\infty} C(\chi_{U^n}) = 0$
	for Outcome \ref{itm:outcome_accumulation_stationary}.
	We do so in two steps. We first prove
	$\liminf_{n\to\infty} C(\chi_{U^n}) = 0$ and
	then improve upon this finding to
	$\lim_{n\to\infty} C(\chi_{U^n}) = 0$.
	
	\emph{Step 1:} We prove
	\begin{gather}\label{eq:step1}
	\liminf_{n\to\infty} C(\chi_{U^n}) = 0.
	\end{gather}
	To this end, we consider the subsequence of accepted 
	iterates (successful steps) $(n_k)_k$ of 
	\cref{alg:binary_tr_steepest_descent}.
	Using the fundamental theorem of calculus and the notation
	$d^n \coloneqq \chi_{U^n\symdif D^n} - \chi_{U^n}$ for $n \in \N$, we may rewrite the decrease in the objective as
	\[
	J(\chi_{U^n}) - J(\chi_{U^n\symdif D^n})= -(\nabla J(\chi_{U^n}),d^n)_{L^2}
	- \int_0^1 (\nabla J(\chi_{U^n} + t d^n) - \nabla J(\chi_{U^n}), d^n)_{L^2} \dd t.
	\]
	We observe that $d^n(s) \in \{-1,0,1\}$ for a.a.\ $s\in\Omega$ implies
	$\sqrt{\|d^n\|_{L^1}} = \|d^n\|_{L^2}$. The Lipschitz continuity of
	$\nabla J : L^1(\Omega) \to L^2(\Omega)$ (\cref{ass:tr_convergence}, \ref{itm:Jprime_Lipschitz}) with Lipschitz constant $L > 0$ implies
	\[ J(\chi_{U^n}) - J(\chi_{U^n\symdif D^n}) \ge 
	-(\nabla J(\chi_{U^n}),d^n)_{L^2}
	- \frac{L}{2} \|d^n\|_{L^1} \sqrt{\|d^n\|_{L^1}}.
	\]
	As in the proof of \cref{lem:criticality_coercion_implies_tr_contraction}, we observe that
	the estimate
	\[ -(\nabla J(\chi_{U^n}), d^n) \ge C(\chi_{U^n})\frac{\Delta^n}{2\lambda(\Omega)} \]
	holds for the steps $d^n$. Inserting this estimate yields
	\[ J(\chi_{U^n}) - J(\chi_{U^n\symdif D^n}) \ge 
	-\sigma_1(\nabla J(\chi_{U^n}),d^n)_{L^2}
	\underbrace{-(1-\sigma_1)\frac{C(\chi_{U^n})\Delta^n}{2\lambda(\Omega)}
		- \frac{L}{2} (\Delta^n)^{\frac{3}{2}}}_{\coloneqq r^n }.
	\]
	We show \eqref{eq:step1} by contradiction.
	If $\liminf_{n\to\infty} C(\chi_{U^n}) > 0$, then \cref{lem:criticality_coercion_implies_tr_contraction}
	implies that $\Delta^n \to 0$ and $r^n \ge 0$ holds for all $n \ge n_2$
	for some $n_2 \in\N$. But $r^n \ge 0$ implies that the acceptance criterion
	in \cref{alg:binary_tr_steepest_descent}, \ref{ln:sufficient_decrease_condition}
	is satisfied for all $n\ge n_2$ and the trust-region radius is not decreased further
	from iteration $n_2$ on. This contradicts $\Delta^n \to 0$ and we thus obtain \eqref{eq:step1}.
	
	\emph{Step 2:} We prove
	\begin{gather}\label{eq:step2}
	\lim_{n\to\infty} C(\chi_{U^n}) = 0.
	\end{gather}
	To this end,
	we follow the proof strategy of \cite[Theorem\,6]{toint1997non}.
	We say that $n \in \N$ is a successful iteration of \cref{alg:binary_tr_steepest_descent} if the acceptance
	test in Line \ref{ln:sufficient_decrease_condition} is successful. Let $\mathcal{S} \subset \N$ denote
	the set of successful iterations of \cref{alg:binary_tr_steepest_descent}.
	We observe that any successful iteration satisfies
	\[ J(\chi_{U^{n}}) - J(\chi_{U^{n+1}}) \ge \frac{1}{2\lambda(\Omega)}C(\chi_{U^n})\Delta^n \]
	because of the properties of the subroutine \texttt{FindStep} (see the proof of
	\cref{lem:criticality_coercion_implies_tr_contraction}). 
	
	This implies
	\[ J(\chi_{U^{0}}) - J(\chi_{U^{n+1}})
	\ge\frac{1}{2\lambda(\Omega)}\sum_{\ell=0, \ell \in S}^n C(\chi_{U^{\ell}})\Delta^\ell
	\]
	for all $n \in \N$. We seek for a contradiction to the claim and assume that there exists a subsequence
	$(n_k)_k \subset S$ such that
	\begin{gather}\label{eq:thm27_contradictory_ass}
	C(\chi_{U^{n_k}}) \ge 2 \varepsilon > 0
	\end{gather}
	for some $\varepsilon > 0$.
	Let $K \coloneqq \{ n \in S\,|\, C(\chi_{U^n}) \ge \varepsilon \}$. It follows that
	\[ J(\chi_{U^{0}}) - J(\chi_{U^{n+1}}) \ge \frac{1}{2\lambda(\Omega)}
	\sum_{\ell=0, \ell \in K}^n C(\chi_{U^{\ell}})\Delta^\ell 
	\ge \frac{1}{2\lambda(\Omega)} \varepsilon \sum_{\ell=0, \ell \in K}^n\Delta^\ell \]
	for all $n \in \N$. Let $n_0 \in \N$. Then we obtain for all $n \ge n_0$ that
	\[ \sum_{\ell=n_0, \ell \in K}^n \Delta^\ell
	\le \frac{2\lambda(\Omega)}{\varepsilon}\bigl(J(\chi_{U^{0}}) - J(\chi_{U^{n+1}}) \bigr) 
	\le\frac{2\lambda(\Omega)}{\varepsilon}\bigl(J(\chi_{U^{0}}) - \min \eqref{eq:r} \bigr) < \infty,
	\]
	which implies that $\sum_{\ell=n_0,\ell \in K}^\infty \Delta^\ell < \kappa < \infty$ for some $\kappa > 0$
	for all $n_0 \in \N$.
	From \eqref{eq:step1} it follows for all
	$k \in \N$ that there exists a smallest
	$\ell(k) > n_k$ with $\ell(k) \in \mathcal{S} \setminus K$.
	We obtain
	\begin{align*}
	\|\chi_{U^{\ell(k)}} - \chi_{U^{n_k}}\|_{L^1}
	&\le \sum_{j=n_k,j \in \mathcal{S}}^{\ell(k) - 1} \|\chi_{U^{j+1}} - \chi_{U^{j}}\|_{L^1}
	\le \sum_{j=n_k,j \in \mathcal{S}}^{\ell(k) - 1} \Delta^j
	= \sum_{j=n_k,j \in K}^{\ell(k) - 1} \Delta^j
	\le \sum_{j=n_k,j \in K}^{\infty} \Delta^j < \kappa.
	\end{align*}
	This implies that for $k \to \infty$ we obtain that
	\[ \|\chi_{U^{\ell(k)}} - \chi_{U^{n_k}}\|_{L^1} \to 0. \]
	By virtue of Fatou's lemma and the fact that every sequence that converges in $L^1$ has a pointwise a.e.\
	convergent subsequence, we obtain
	\[ \|\chi_{U^{\ell(k)}} - \chi_{U^{n_k}}\|_{L^2} \to 0 \]
	for a subsequence of $(n_k)_k$, which we denote with the same symbol for ease of notation.
	We conclude that
	\[ |C(\chi_{U^{\ell(k)}}) - C(\chi_{U^{n_k}})| \to 0~ (\text{for }k\ \to \infty),
	\]
	which violates our assumption \eqref{eq:thm27_contradictory_ass}. Hence, \eqref{eq:step2} follows.
	
	For the remainder of the proof, we restrict ourselves to 
	a weakly convergent subsequence of $(\chi_{U^n})_n$,
	for ease of notation denoted by
	the same symbol, which satisfies
	$C(\chi_{U^{n}}) \to 0$ and $\chi_{U^n} \rightharpoonup f$ in $L^2(\Omega)$.
	It remains to show that $C(f) = 0$.
	if Assumptions \ref{ass:tr_convergence},
	\ref{itm:Jdifferentiable}
	and \ref{itm:Jprime_cc} hold.
	The criticality measure $C$ is weakly lower semi-continuous
	under \cref{ass:tr_convergence}, \ref{itm:Jprime_cc}, see \cite[Lemma 4.1]{dunn1980convergence},
	so that 
	\[ 0 \le C(f)
	\le \liminf_{n\to\infty} C(\chi_{U^n}) = 0. 
	\]
\end{proof}

\section{Computational Experiments}\label{sec:compu-valid}

We carry out our experiments on an instance of \eqref{eq:p} that
satisfies \cref{ass:tr_convergence}. The instance is described in 
\S\ref{sec:example}, and our computational setup is described in 
\S\ref{sec:setup}. Then validation experiments and their results for
the presented theory are presented in \S\ref{sec:validation}.
Motivated by observations in \S\ref{sec:validation}, we explore
the effects of a hybridization of BTR and CIA in
\S\ref{sec:hybridization}.

\subsection{Example Problem}\label{sec:example}

We consider the case $d = 2$ and the domain $\Omega = (0,2)^2 \subset \R^2$.
For \eqref{eq:p} we choose $J(x) \coloneqq j(S(x))$,
where $j$ is a so-called tracking-type objective, specifically
$j(y) = .5\|y - y_d\|_{L^2}^2$ for a given
\[ y_d(s) = \frac{1}{4} \sin(3 (s_1 - 1)(s_2 - 1))^2 (|s_1 - 1| + |s_2 - 1|) \]
for $s \in \Omega$. $S$ is the solution operator of the linear elliptic boundary 
value problem
\begin{equation}\label{eq:poisson-bvp}
-\varepsilon \Delta y + y = x,\quad y|_{\partial\Omega} = 0
\end{equation}
for a given control input $x \in L^1(\Omega)$ and the choice $\varepsilon = 10^{-2}$.
This yields the following instance of \eqref{eq:p}:
\begin{equation}\label{eq:poisson-tracking}
\inf_x\quad \frac{1}{2}\|y - y_d\|_{L^2}^2
\quad\text{s.t.}\quad y = S(x) \text{ and } x(s) \in \{0,1\} \text{ for a.a.\ } s \in \Omega.
\end{equation}
For a Poisson problem with  right-hand side in $L^1(\Omega)$,
the weak solution $y$ is an element of the Sobolev space of $q$-integrable
functions with a $q$-integrable distributional derivative that vanish at 
the boundary, $W^{1,q}_0(\Omega)$, for bounded Lipschitz domains 
$\Omega$, where we have the estimate $\|u\|_{W^{1,q}_0(\Omega)} \le c\|x\|_{L^1(\Omega)}$ for some $c > 0$ if $q < 2$.
A proof of this result (for more general elliptic operators) and
more general right-hand sides can, for example, be found in 
\cite[Theorem 1]{boccardo1989non}; and for the case of mixed
boundary condition a proof can be found in \cite{haller2009holder}
(note that the required Gr\"{o}ger regularity of the boundary therein 
reduces to requiring a strong Lipschitz condition if only a Dirichlet 
boundary condition is present).

Combining these considerations with the chain rule for Banach spaces and the
Riesz representation theorem implies that
Assumptions \ref{ass:tr_convergence},
\ref{itm:Jboundedbelow} and \ref{itm:Jprime_Lipschitz}
are satisfied for this example.
To see that \cref{ass:tr_convergence} \ref{itm:Jprime_cc} is also satisfied, we
consider the compact embedding
$W^{1,q}_0(\Omega) \hookrightarrow W^{1,1}(\Omega) \hookrightarrow^c L^2(\Omega)$
(for $d = 2$), where the compactness is due to the second embedding.
Because $S'(\bar{x})$ is linear and bounded for
$\bar{x} \in L^1(\Omega)$, it maps weakly convergent sequences to weakly convergent sequences
in $W^{1,q}_0(\Omega)$ and, by compactness, to norm convergent sequences in $L^2(\Omega)$, which implies that $S'$ and in turn $\nabla J$ are weak-norm
continuous.

\subsection{Setup}\label{sec:setup}

We solve the boundary value problem \eqref{eq:poisson-bvp} numerically 
using a finite element method on a conforming uniform triangle mesh
subdividing the  domain $\Omega$. Specifically, the domain is
partitioned into $256 \times 256$ square cells that are split into
$4$ triangles each. In all experiments, the solution $y = S(x)$ of 
\eqref{eq:poisson-bvp} is computed in the space of cellwise affine
and globally continuous functions on this mesh.
All experiments are carried out on a laptop computer with Intel(R) Core(TM) 
i9-10885H CPU (2.40\,GHz) and 32 GB RAM.

\paragraph{Implementation of CIA
	(\cref{alg:cia})}

We compute a 
solution of the continuous relaxation of \eqref{eq:poisson-tracking}
(the first step of CIA)
by replacing the constraint $x(s) \in \{0,1\}$ by $x(s) \in [0,1]$ and
optimizing the control function on the aforementioned triangle mesh 
in the space of cellwise constant discontinuous functions with a 
quasi-Newton method. 

For the \texttt{Round} procedure, the second 
step of CIA, we consider three different choices:
multidimensional sum-up rounding (SUR) \cite{manns2020multidimensional},
the combinatorial optimization-based rounding (COR)
proposed in \cite{jung2015lagrangian}, see also
\cite[\S\,2.4.1]{bestehorn2022combinatorial},
and a primal heuristic for switching cost aware rounding (SHG)
\cite{bestehorn2019switching,bestehorn2021mixed}.
The choices for the \texttt{Round}
procedure compute approximating controls. In our experiments, the
computed controls are cellwise constant functions
on cells of the grid
of $256\times 256$ squares. The three choices
for the \texttt{Round}
procedure are explained in \S\ref{sec:sur-appendix}.
In order to ensure the approximation property 
\eqref{eq:round_weakstar}, we order
the grid cells of the discretization along a Hilbert 
curve as in 
\cite{manns2020multidimensional,kirches2020compactness},
which yields an order-conserving domain dissection
as mentioned in \S\ref{sec:cia},
see also \cref{dfn:order-conserving-domain-dissection}.

While SUR and COR either minimize
a certain approximation error or abide by an upper bound on that error,
SHG minimizes the length of the interface between the level sets 
for the values zero and one in our 
setting while abiding by the approximation error bound.
The integer programming formulation of SHG is computationally intractible
in our setting (for our grid size), which is why
we resort to a suboptimal heuristic, see 
also \S\ref{sec:sur-appendix}.

\paragraph{Implementation of BTR
	(\cref{alg:binary_tr_steepest_descent})} We compute all iterates on the fixed 
uniform mesh of  $256 \times 256$ squares. The uniformity, namely, the
fact that all cells have the same volume, has the advantage that
a discretized variant of the \texttt{FindStep} method can be implemented
efficiently. Specifically it is a Knapsack problem with all
weights being equal to one, which is therefore not NP-hard.
In particular, if $x = \chi_A$ for some $A \in \B$ is the current
iterate, the discretized trust-region subproblem \eqref{eq:tr}
can be solved as follows. We compute the average value of the
function $g = \nabla J(\chi_A)(\chi_{A^c} - \chi_A)$ on
each grid cell. Then we sort the cell averages of $g$ in ascending
order and pick the cells with negative average values in a greedy 
fashion until the current trust region is filled. The picked cells 
constitute a difference set $D$ so that the computed step $d$ is
$d = \chi_{A \symdif D} - \chi_A$. Our discretized implementation
of BTR terminates when the trust-region radius contracts below the
volume of one grid cell.

We note that while we have carried out our experiments on a fixed fine 
mesh, an alternative approach is to adaptively refine the mesh where 
required within the \texttt{FindStep} 
subroutine of BTR.

\subsection{Validation}\label{sec:validation}

We apply CIA to \eqref{eq:poisson-tracking}, where we initialize the
solver for the continuous relaxation with the constant zero function, 
$x=0$. As mentioned above, we use SUR, COR, and SHG for the \texttt{Round} subroutine
in \cref{alg:cia} in order to compute the binary-valued approximation of the solution
of the continuous relaxation on the uniform mesh of $256 \times 256$ squares.
COR and SUR produce the same resulting control. SHG also returns the
same control if the approximation constraint in the problem formulation is used
with $\theta = 1$, see \cref{alg:shg}, indicating that the feasible set with this bound
leaves (almost) no room for a reduction of the interface length. We therefore increase
the feasible set by setting $\theta = 10$, which provides a trade-off that relaxes the 
approximation quality and allows to reduce the interface length within the prescribed
approximation quality. The optimality gap is then higher but the interface length decreases
and we obtain a qualitatively different solution.

Then we initialize BTR, which also operates on the
uniform grid of $256 \times 256$ squares, with the output
control of SUR. Because of the near-optimality achieved by CIA and
the fact that BTR does not produce a globally optimal
solution of \eqref{eq:poisson-tracking} for a fixed discretization,
we expect that BTR can close a small part but not much of the
remaining optimality gap between the upper bound given by the objective 
value for the output SUR and the lower bound given by the solution of
the continuous relaxation. 

This expectation is met by our computational results. Specifically,
BTR is able to close a portion of the remaining optimality gap of CIA.
We also start BTR from two further different initializations,
specifically, from $x = 0$ and from a cellwise rounding of the solution of the
continuous relaxation to $\{0,1\}$. All of the objective values are
very close; the optimality gap is always around $10^{-6}$
(at a magnitude of $10^{-3}$) so that BTR also produces a 
near-optimal solution on this grid when it is initialized differently.

Moreover, we observe that the running times of BTR are longer than those of
CIA, which we attribute to the facts that BTR is a pure first-order 
method, has restricted options for its feasible steps available 
(compared with usual gradient-based solvers for the continuous 
relaxation), and requires a sorting operation additionally to each
adjoint solve of each accepted step. The higher running times are
reflected by correspondingly high numbers of iterations of
our implementation of \cref{alg:binary_tr_steepest_descent},
as is typical for first-order methods. The running time of the
continuous relaxation followed by an execution of BTR on
cellwise rounding is moderately lower (about $20\%$) than that of
BTR for initial control zero.

\begin{figure}
	\centering
	\subfloat[Rel.]{
		\includegraphics[width=3.75cm]{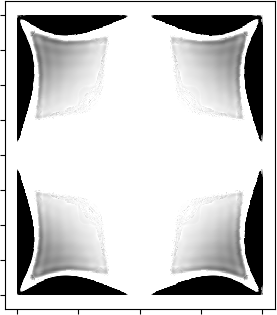}}\hspace{0.5cm}
	\subfloat[CIA (SUR) / (COR)]{
		\includegraphics[width=3.75cm]{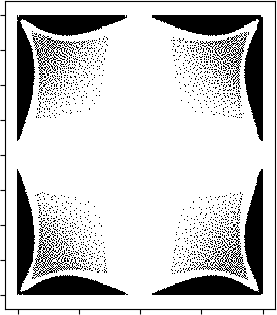}}\hspace{0.5cm}
	\subfloat[CIA (SHG)]{
		\includegraphics[width=3.75cm]{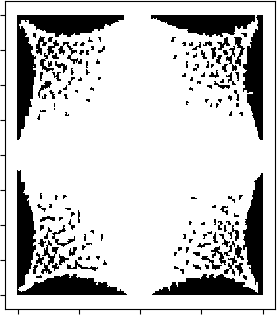}}\hspace{0.5cm}    
	\subfloat[CIA (SUR) + BTR]{
		\includegraphics[width=3.75cm]{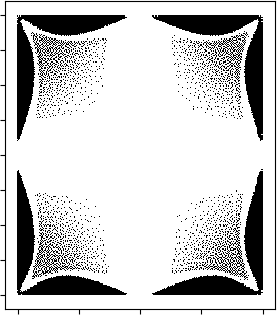}}\hspace{0.5cm}
	\subfloat[Rel.\ + BTR]{
		\includegraphics[width=3.75cm]{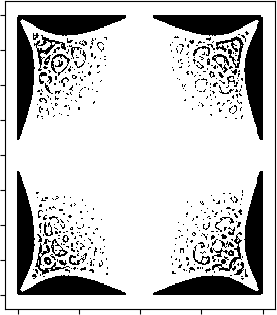}}\hspace{0.5cm}
	\subfloat[BTR]{
		\includegraphics[width=3.75cm]{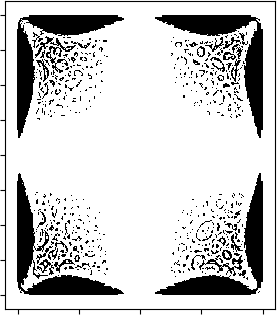}}
	\caption{Visualization of the resulting control functions
		for \emph{Rel.}
		(= continuous relaxation of \eqref{eq:poisson-bvp}),
		\emph{CIA (SUR) / (COR)} (= continuous relaxation and SUR / COR),
		\emph{CIA (SHG)} (= continuous relaxation and SHG),    
		\emph{CIA (SUR) / (COR) + BTR}, \emph{Rel.\ + BTR}
		(= continuous relaxation and BTR started from a cellwise rounding),
		and \emph{BTR} (= BTR started from zero), where
		the value one is colored black and the value zero
		is colored white. For \emph{Rel.}, the intermediate
		values in $[0,1]$ are depicted in grayscale.}
	\label{fig:validation}
\end{figure}

To give a qualitative and visual impression of the
results, we provide the six computed controls in this
experiment in \cref{fig:validation}.
A visual inspection of the results shows that starting 
BTR from a cellwise rounded solution or zero seems to 
have a regularizing effect
on the resulting microstructure. We
compute the length of the interface between
the level sets for the values zero and the objective
value one for all
of the computed controls. We obtain that the none
of the controls computed with
CIA with SUR or COR as the \texttt{Round} procedure,
CIA with SHG as the \texttt{Round} procedure,
and BTR started from cellwise rounding
dominates another in terms of low objective
or low interface length value. Specifically,
CIA with SUR or COR produces the lowest objective
value but the highest interface length.
The interface lengths for the controls
with BTR started from zero or cellwise rounding
are significantly lower (about 50\,\%) while
the objective value increases slightly.
The interface length decreases further (about 30\,\%)
for CIA with SHG for $\theta = 10$
but the increase in the objective
is also higher.

While the behavior of the
objective values follows our analysis, we cannot
explain our consistent observation that some spatial
coherence is maintained during the optimization
with BTR.
In order to do so, we need to better
understand the possible trade-off between objective
and interface length, which also motivated us
to introduce SHG and report its results, which is
impaired by its current computational intractability 
and us resorting to a heuristic solution. We leave 
further considerations in this direction to future 
research.

These findings also lead to the question 
whether one can obtain the
regularization effect in a hybridized method, where BTR is 
initialized with an approximation of SUR that is computed on a 
coarser control mesh. This is investigated in 
\S\ref{sec:hybridization}.
We report the objective values, optimality
gaps, running times,  number of iterations required by
our implementation of BTR, and  interface lengths
between the level sets in \cref{tbl:validation}.

\begin{table}
	\setlength\tabcolsep{5pt}
	\renewcommand{\arraystretch}{.95}
	\centering
	\caption{Objective values, remaining optimality gaps, 
		BTR iteration numbers, running times, and interface
		lengths obtained for the cases \emph{Rel.}
		(= continuous relaxation of \eqref{eq:poisson-bvp}),
		\emph{CIA (COR)} (= continuous relaxation and COR),
		\emph{CIA (SUR)} (= continuous relaxation and SUR),    
		\emph{CIA (SHG)} (= continuous relaxation and SHG),        
		\emph{CIA + BTR} (= continuous relaxation and BTR started from SUR),
		\emph{Rel.\ + BTR}
		(= continuous relaxation and BTR started from a cellwise rounding),
		and \emph{BTR} (= BTR started from zero).}\label{tbl:validation}
	\begin{tabular}{r|lllllll}
		\toprule
		& Rel. 
		& CIA (COR)
		& CIA (SHG)
		& CIA (SUR)
		& CIA + BTR & Rel.\ + BTR & BTR  \\
		\midrule                
		Obj.\ [$10^{-3}$]
		& $4.0798$
		& $4.0808$
		& $4.1200$
		& $4.0808$
		& $4.0807$
		& $4.0837$
		& $4.0862$ \\
		Opt.\ gap [$10^{-6}$]  
		& 0
		& 1.06
		& 40.20
		& 1.06
		& 0.89
		& 3.96
		& 6.41 \\
		Time [s]        & 821 & 828 & 829 & 825 & 921 & 3928 & 4906 \\
		BTR iterations  & n/a & n/a & n/a & n/a & 51 & 1575 & 2481 \\
		Interface length & n/a & 117.0 & 48.1 & 117.0 & 116.9 & 66.4 & 74.2 \\
		\bottomrule
	\end{tabular}
\end{table}

We also evaluate how our implementation of BTR behaves with respect to
mesh refinement when initialized with zero und running it until the trust-region radius
contracts. The optimality gap decreases for finer meshes as finer
microstructures can be computed in order to more closely approximate a minimizer of the continuous
relaxation. However, the number of iterations grows with a factor of approximately three when the mesh
size is halved (the number of cells increases by a factor of four).
In particular, our implementation of the algorithm is 
clearly not mesh-independent.
The interface length also increases when finer meshes are chosen, which is 
consistent with the fact that this quantity tends to infinity if a function
with values in $(0,1)$ on a set of strictly positive measure
(in this case the solution of the relaxation)
is approximated weakly-$^*$ in $L^\infty$ by binary functions.
We provide the corresponding data in \cref{tbl:mesh_refinement}.

\begin{table}
	\setlength\tabcolsep{5pt}
	\renewcommand{\arraystretch}{.95}
	\centering
	\caption{Remaining optimality gaps, BTR iteration numbers, and interface
		lengths for BTR started from zero for different mesh sizes.}\label{tbl:mesh_refinement}
	\begin{tabular}{r|lll}
		\toprule
		Mesh & Optimality gap $[10^{-6}]$ & BTR iterations & Interface length  \\
		\midrule                
		$32 \times 32$ & 1200.79 & 116 & 21.8 \\
		$64 \times 64$ & 261.14 & 356 & 34.2 \\
		$128 \times 128$ & 38.02 & 938 & 50.6 \\
		$256 \times 256$ & 6.41 & 2481 & 74.2 \\
		\bottomrule
	\end{tabular}
\end{table}

\subsection{Hybridization of SUR and BTR}\label{sec:hybridization}

We explore the regularization effect that we have observed
above by executing a hybridized method by initializing BTR with 
controls that are computed by CIA, where the second step
is computed with SUR.

SUR operates on a mesh, and its approximation quality depends
on the mesh size of this grid
\cite{kirches2020compactness,manns2020multidimensional}.
We use a sequence of uniformly refined grids from $8\times 8$
grid cells to $256 \times 256$ grid cells
and start BTR, which itself operates on the  $256 \times 256$ grid
of squares, with the resulting controls. 

We assess the running times and iterations of BTR, the remaining 
optimality gaps, and the interface lengths between the 
level sets of the two control realizations one and zero.

The remaining optimality gaps are of an order of magnitude
of $10^{-6}$, where  the ones achieved for the initializations
with the SUR solutions for the $128\times 128$ and $256\times 256$
grids are slightly but noticeably smaller.
This can be attributed to the
fact that SUR provides already small optimality gaps, which
are then further reduced by BTR.

Moreover, the running times and iterations of BTR do not
differ much from the $8\times 8$ to $64\times 64$ initializations
(all above 3000\,s) and then drop in two large steps to
1074\,s and 96\,s for the two finest grids,
again capitalizing on the fact that SUR already provides a
near-optimal solution.

The interface lengths of the resulting level sets of the controls
obtained for the $8\times 8$ to $128\times 128$ initializations
are between 62.5 and 70.2, approximately 50\,\% smaller
than the interface length obtained for the finest grid.
This indicates that one may find a sensible trade-off,
where the regularization effect of BTR is still pronounced
and BTR can be meaningfully accelerated by means of CIA,
in our case by executing SUR in CIA on the
$128\times 128$ initialization.

\begin{figure}
	\centering
	\subfloat[SUR ($8\times 8$)]{
		\includegraphics[width=3.75cm]{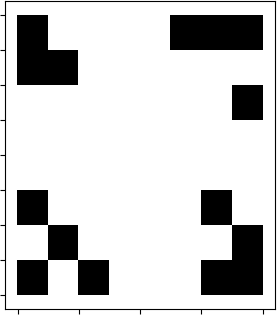}}\hspace{0.5cm}
	\subfloat[SUR ($128\times 128$)]{
		\includegraphics[width=3.75cm]{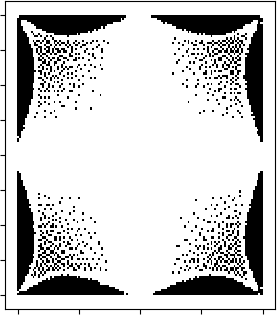}}\\
	\subfloat[BTR from SUR ($8\times 8$)]{
		\includegraphics[width=3.75cm]{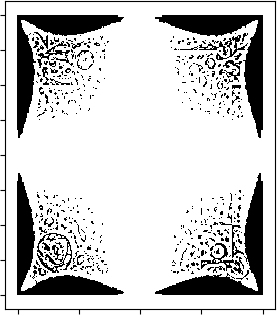}}\hspace{0.5cm}
	\subfloat[BTR from SUR ($128\times 128$)]{
		\includegraphics[width=3.75cm]{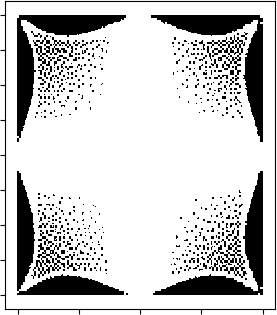}}
	\caption{Visualization of the resulting control functions
		for SUR on uniform $8\times 8$ and $128\times 128$ square
		grids (top row) and BTR initialized with them
		(bottom row).}
	\label{fig:hybridization}
\end{figure}

We have recorded the obtained results in \cref{tbl:hybridization}.
To provide a visual impression again, we
contrast the computed controls for SUR with the resulting
ones after initializing BTR with them for the
$8\times 8$ and $128\times 128$ grids in \cref{fig:hybridization}.
\begin{table}[h]
	\setlength\tabcolsep{5pt}
	\renewcommand{\arraystretch}{.95}
	\centering
	\caption{Remaining optimality gaps,
		running times, BTR iterations, and interface lengths
		for executing BTR on solutions of CIA, where
		the second step is computed by means of SUR,
		for refined grids used for SUR.}\label{tbl:hybridization}
	\begin{tabular}{r|lllllll}
		\toprule
		SUR Mesh
		& $8\times 8$
		& $16\times 16$
		& $32\times 32$
		& $64\times 64$
		& $128\times 128$ 
		& $256 \times 256$ \\
		\midrule                
		Optimality gap [$10^{-6}$]  
		& 4.99
		& 6.66
		& 4.68
		& 5.30
		& 2.91
		& 0.89 \\
		Time (only BTR) [s]
		& 4376 
		& 3929
		& 4049
		& 3167
		& 1074 
		& 96\\
		BTR iterations  & 2242
		& 2023
		& 2069
		& 1619
		& 552
		& 51 \\
		Interface length 
		& 70.2 
		& 67.9
		& 67.4
		& 62.5 
		& 67.0
		& 116.9 \\
		\bottomrule
	\end{tabular}
\end{table}

\section*{Acknowledgments}

We thank two anonymous referees for providing helpful
feedback on the manuscript. We thank Peter Bella and 
Christian Meyer (both TU Dortmund University) for helpful 
discussions on the topic. This work was supported by the U.S.~Department of Energy, Office of Science, Office of Advanced Scientific Computing Research, Scientific Discovery through the
Advanced Computing (SciDAC) Program through the FASTMath Institute under Contract No. DE-AC02-06CH11357 and by the German Research Foundation under GRK 2297 MathCoRe (project No. 314838170), SPP~1962 (projects No. SA 2016/1-2, KI 1839/1-2), and SPP~2231 (project No. SA 2016/3-1, KI 417/9-1), and by the German Federal Ministry of Education and Research within the program ``Mathematics for Innovations'', project ``Power to Chemicals.'' This work is also part of the Research Initiative ``SmartProSys: Intelligent Process Systems for the Sustainable Production of Chemicals'' funded by the Ministry for Science, Energy, Climate Protection and the Environment of the German State of Saxony-Anhalt.

\appendix

\section{Auxiliary Results}\label{sec:appendix}
This appendix provides additional information on the \texttt{Round} subroutine in the CIA method, \cref{alg:cia}, and establishes 
relationships between set-based and characteristic function.
We also prove the equivalence between \eqref{eq:g_subproblem}
and \eqref{eq:tr} and the existence of minimizers for them.

\subsection{\texttt{Round} in \cref{alg:cia}}\label{sec:sur-appendix}
\Cref{alg:cia} employs a \texttt{Round} subroutine, which takes grids as inputs. In order to establish the 
approximation property \eqref{eq:round_weakstar}, see also \cite{manns2020multidimensional,kirches2020compactness},
it is assumed that the sequence of grids is a so-called \emph{order-conserving domain
	dissection}, which is formally defined below.
\begin{definition}[Order-conserving domain dissection, Definition 4.3 in \cite{manns2020multidimensional}]\label{dfn:order-conserving-domain-dissection}
	Let $\Omega \subset \R^d$. Then we call a sequence
	$\left(\{S_1^n,\ldots,S^n_{N^n}\}\right)_n \subset 2^{\mathcal{B}(\Omega)}$ 
	an \emph{order-conserving domain dissection} if
	\begin{enumerate}
		\item $\{S_1^n,\ldots,S^n_{N^n}\}$ is a finite partition of $\Omega$ for all $n \in \N$,
		\item $\max\{ \lambda(S_i^n) \,|\, i \in \{1,\ldots,N^n\}\} \to 0$,
		\item for all $n \in\N$ for all $i \in \{1,\ldots,N^{n-1}\}$ there exists
		$1 \le j < k \le N^n$ such that $\bigcup_{\ell=j}^k S^n_\ell = S^{n-1}_i$, and
		\item the cells $S^n_j$ shrink regularly (there exists $C > 0$ such that
		for each $S_j^n$ there exists a ball $B_j^n$ with $S_j^n \subset B_j^n$
		and $\lambda(S_j^n) \ge C \lambda(B_j^n)$).
	\end{enumerate}
\end{definition}
We briefly introduce the three choices for \texttt{Round} in \cref{alg:cia}, Line \ref{ln:round} that are used
in our computational experiments in \S\ref{sec:compu-valid}  below.
\paragraph{Sum-up rounding (SUR)}
The algorithm  is stated as \cref{alg:sur} and works as follows.
\begin{algorithm}
	\caption{Sum-Up Rounding (multidimensional variant) \cite{manns2020multidimensional}}\label{alg:sur}
	\textbf{Input:} Ordered grid cells $S_1,\ldots,S_N \subset \Omega$ that partition $\Omega$ .
	
	\textbf{Input:} Function $\alpha \in L^1(\Omega,\R^m)$ with averages $a_{k,i}$
	such that $a_{k,i} = \frac{1}{\lambda(S_k)} \int_{S_k} \alpha_i(s)\dd s$
	for all $i\in \{1,\ldots,m\}$, and $\sum_{i=1}^m a_{k,i} = 1$ for all $k \in \{1,\ldots,N\}$.
	\begin{algorithmic}[1]
		\State $\phi_0 \coloneqq 0_{\R^m}$
		\For{k = 1,\ldots,N}
		\State $\gamma_k \gets \phi_{k-1} + a_{k} \lambda(S_k)$\quad\quad\quad\quad\ {\footnotesize ($a_k$ is short for the vector $(a_{k,1},\ldots,a_{k,m})^T$)}
		\State $w_{k,i} \gets \left\{ \begin{aligned}
		1 & : && i \in \argmax \{ \gamma_{k,j}\,|\, j \in \{1,\ldots,m\} \},\\
		0 & : && \text{else}
		\end{aligned}\right.$ for all $i \in \{1,\ldots,m\}$
		\State $\phi_k \gets \sum_{i=1}^k (a_k - w_k) \lambda(S_k)$\quad\quad {\footnotesize ($w_k$ is short for the vector $(w_{k,1},\ldots,w_{k,m})^T$)}
		\EndFor
		\State \textbf{return} $\omega \coloneqq \sum_{k=1}^N w_k \chi_{S_k}$
	\end{algorithmic}
\end{algorithm}
It starts from a $[0,1]^m$-valued
function $\alpha$ such that
$\sum_{i=1}^m \alpha_i = 1$ a.e.\ and that is defined on an ordered sequence of grid cells that
partition the domain $\Omega$. The algorithm iterates over the grid cells in the given order; identifies
an entry $i \in \{1,\ldots,m\}$ such that the cumulative difference up to the current grid cell to a rounded function
$\omega$, which is $\{0,1\}^m$-valued, satisfies $\sum_{i=1}^m \omega_i = 1$ a.e., and is defined on the same grid as
$\alpha$, is maximal. Then the algorithm sets $\omega_i$ to one on the current grid cell and the other entries
to zero on the respective grid cell.

\paragraph{Combinatorial optimization-based rounding (COR)}
The algorithm is stated as \cref{alg:cor} and works as follows.
\begin{algorithm}
	\caption{Combinatorial optimization-based rounding (multidimensional variant)}\label{alg:cor}
	\textbf{Input:} Ordered grid cells $S_1,\ldots,S_N \subset \Omega$ that partition $\Omega$ .
	
	\textbf{Input:} Function $\alpha \in L^1(\Omega,\R^m)$ with averages $a_{k,i}$
	such that $a_{k,i} = \frac{1}{\lambda(S_k)} \int_{S_k} \alpha_i(s)\dd s$
	for all $i\in \{1,\ldots,m\}$, and $\sum_{i=1}^m a_{k,i} = 1$ for all $k \in \{1,\ldots,N\}$.
	\begin{algorithmic}[1]
		\State compute $w$ as minimizer of
		
		$\underset{w,\eta}\argmin\ \eta\quad\text{s.t.}\quad
		\left\{\begin{aligned}
		& w \in \{0,1\}^{N\times m}\\
		&\sum_{i=1}^m w_{k,i} = 1 \text{ for all } k \in \{1,\ldots,N\}\\
		& \left|\sum_{j=1}^k (a_{j,i} - w_{j,i}) \lambda(S_j)\right| \le \eta
		\text{ for all } k \in \{1,\ldots,N\} 
		\text{ and all } i \in \{1,\ldots,m\}
		\end{aligned}\right.
		$
		\State \textbf{return} $\omega \coloneqq \sum_{k=1}^N w_k \chi_{S_k}$
	\end{algorithmic}
\end{algorithm}
It starts from a $[0,1]^m$-valued
function $\alpha$ such that
$\sum_{i=1}^m \alpha_i = 1$ a.e.\ and that is defined on an ordered sequence of grid cells that
partition the domain $\Omega$. The algorithm computes a $\{0,1\}^m$-valued function that satisfies 
$\sum_{i=1}^m \omega_i = 1$ a.e.\ that is piecewise constant on this grid. The function is computed
such that it minimizes the maximum of the modulus of the cumulative difference (integration)
to $\alpha$ from the first to the $k$-th cell over $k \in \{1,\ldots,N\}$.
\Cref{alg:cor} can be implemented very efficiently using the shortest path approach described
in \cite{bestehorn2022combinatorial,bestehorn2021mixed} on multidimensional
domains because the costs are sequence independent in the sense of \cite[\S\,2.6.1]{bestehorn2022combinatorial}.
We use the open-source software \texttt{scarp\_solver}\footnote{\label{ftn:gh}
	Accessed on \url{https://github.com/chrhansk/SCARP} on 02/15/2022.}
with the option \texttt{-{}-sur\_costs} in order to solve SHG in this work.

\paragraph{Primal heuristic for switching cost aware rounding (SHG)}
Switching cost aware rounding \cite{bestehorn2021mixed,bestehorn2021switching}
is stated as \cref{alg:shg} and works as follows.
\begin{algorithm}
	\caption{Switching cost aware rounding (multidimensional variant)}\label{alg:shg}
	\textbf{Input:} Ordered grid cells $S_1,\ldots,S_N \subset \Omega$ that partition $\Omega$.
	
	\textbf{Input:} Function $\alpha \in L^1(\Omega,\R^m)$ with averages $a_{k,i}$
	such that $a_{k,i} = \frac{1}{\lambda(S_k)} \int_{S_k} \alpha_i(s)\dd s$
	for all $i\in \{1,\ldots,m\}$, and $\sum_{i=1}^m a_{k,i} = 1$ for all $k \in \{1,\ldots,N\}$.
	
	\textbf{Input:} Trade-off parameter $\theta \ge 1$.
	\begin{algorithmic}[1]
		\State compute $w$ as minimizer of
		
		$\underset{w}\argmin\ \frac{1}{m} \TV\left(\sum_{k=1}^N w_k \chi_{S_k}\right) \quad\text{s.t.}\quad
		\left\{\begin{aligned}
		& w \in \{0,1\}^{N\times m}\\
		&\sum_{i=1}^m w_{k,i} = 1 \text{ for all } k \in \{1,\ldots,N\}\\
		& \left|\sum_{j=1}^k (a_{j,i} - w_{j,i}) \lambda(S_j)\right| \le \theta \sum_{i=2}^m\frac{1}{i} \max_{\ell}\lambda(S_\ell)\\
		&\quad\quad\text{for all } k \in \{1,\ldots,N\}  \text{ and all } i \in \{1,\ldots,m\}
		\end{aligned}\right.
		$
		\State \textbf{return} $\omega \coloneqq \sum_{k=1}^N w_k \chi_{S_k}$
	\end{algorithmic}
\end{algorithm}
It starts from a $[0,1]^m$-valued
function $\alpha$ such that
$\sum_{i=1}^m \alpha_i = 1$ a.e.\ and that is defined on an ordered sequence of grid cells that
partition the domain $\Omega$. The algorithm computes a $\{0,1\}^m$-valued function that satisfies 
$\sum_{i=1}^m \omega_i = 1$ a.e.\ that is piecewise constant on this grid. The function is computed
such that it minimizes the total variation of $w$ while constraining the modulus of the cumulative difference
(integration) to $\alpha$ from the first to the $k$-th cell over $k \in \{1,\ldots,N\}$
by $\theta \sum_{i=2}^m\frac{1}{i} \max_{\ell}\lambda(S_\ell)$.
$\TV(w)$ denotes the total variation of $w$ in \cref{alg:shg}.
The constant $\sum_{i=2}^m\frac{1}{i} \max_{\ell}\lambda(S_\ell)$ is guaranteed by \cref{alg:sur},
see the analysis in \cite{kirches2020approximation,manns2021approximation}, so that a feasible point always
exists for $\theta = 1$. The feasible set may be increased by choosing $\theta > 1$ in order to leave
room for a better objective while allowing for a larger approximation error.
After discretization, the optimization problem in \cref{alg:shg} becomes an integer linear
program. While a shortest path reformulation and efficient solution algorithms exist
for the case that $\Omega \subset \R$, it is not known if it can be solved efficiently
for $\Omega \subset \R^d$, $d \ge 2$. By considering only the jumps between subsequent
cells along the ordering of the grid cells in the minimization, we can obtain a suboptimal feasible point
using the shortest path approach described in \cite{bestehorn2022combinatorial,bestehorn2021mixed}.
This is introduced as the heuristic \texttt{SCARP\_HG} in \cite{bestehorn2021switching}.
We use the open-source software \texttt{scarp\_solver}\footnotemark[\getrefnumber{ftn:gh}] with 
the option \texttt{-{}-scale $\theta$} in order to solve \texttt{SCARP\_HG} / (SHG).

\subsection{Equivalence of \eqref{eq:tr}
	to \eqref{eq:g_subproblem} and Existence of Minimizer}

\begin{proposition}\label{prp:tr_equivalence}
	Let $g = \nabla J(\chi_U)(\chi_{U^c} - \chi_U)$.
	A set $D \in \B$ satisfies $\lambda(D) \le \Delta$
	if and only if $d = \chi_{U \symdif D} - \chi_{U}$
	is feasible for \eqref{eq:tr} and
	the corresponding objective values for 
	\eqref{eq:g_subproblem} and \eqref{eq:tr}
	coincide.
	
	If $D \in \B$ does not satisfy $D \subset g^{-1}((-\infty,0])$, then its objective value is greater
	or equal than $D \cap g^{-1}((-\infty,0])$
	so that the optimal objective for \eqref{eq:tr} is not
	altered by the additional feasible points.
\end{proposition}
\begin{proof}
	For every $D \in \B$, $d$ can be computed with the 
	formula above. On the other hand, the constraint
	$\chi_{U}(s) + d(s) \in \{0,1\}$ in \eqref{eq:tr}
	implies that $\chi_{U} + d$ is a characteristic
	function of a measurable set $A$, which in turn can be
	represented as $A = U \symdif D$
	for the set $D = U \symdif A$. Moreover, we have
	$\|d\|_{L^1} = \|\chi_{U \symdif D} - \chi_{U}\|_{L^1} = \lambda(D)$, which shows the equivalence
	of the trust-region constraint.
	
	Moreover, for any feasible $D$ and $d = \chi_{U^n \symdif D} - \chi_{U}$, we distinguish
	the four cases whether $s \in D$ and/or in $s \in U$
	holds and obtain
	\begin{align*}
	\int_D g \dd s
	&= \int_D \nabla J(\chi_{U})(\chi_{U^c} - \chi_{U}) \dd s
	= \int_\Omega \nabla J(\chi_{U})d \dd s
	= (\nabla J(\chi_{U}), d)_{L^2},
	\end{align*}
	which gives the coincidence of the objective
	values.
	
	Finally, let $D \in \B$ be given. Then
	\[
	\int_D g \dd s
	= \int_{D \cap g^{-1}((-\infty,0])} g \dd s 
	+ \int_{D \cap g^{-1}((0,\infty))} g \dd s
	\le \int_{D \cap g^{-1}((-\infty,0])} g \dd s.
	\]
\end{proof}

\begin{proposition}\label{prp:min_tr_subproblem}
	Let $\Omega$ be a bounded domain.
	Let $U \in \B$. Let $\Delta \in [0,\infty]$. Let $\nabla J \in L^1(\Omega)$. Then \eqref{eq:tr} admits a minimizer $\tilde{d}$ with
	$(\nabla J(\chi_U), \tilde{d})_{L^2} \le -\frac{\Delta}{\lambda(\Omega)} C(\chi_U)$.
\end{proposition}
\begin{proof}
	Let $g \coloneqq \nabla J(\chi_U)(\chi_{U^c} - \chi_U)$, and let $D : \R \to \B$ be defined as $D(x) \coloneqq g_U^{-1}((-\infty,x))$ for $x \in \R$. Then $D(x) \subset D(y)$ for all $x \le y$, implying that $\|\chi_{U \Delta D(x)} -\chi_{U}\|_{L^1} = \lambda(D(x))$ is monotone
	in $x$.
	Let $d(x) \coloneqq \chi_{U \Delta D(x)} -\chi_{U}$. Then the $d(x)$ are greedy solution candidates for
	\eqref{eq:tr} with $\lim_{x\to-\infty}\|d(x)\|_{L^1} = 0$. Specifically, $d(0)$ minimizes \eqref{eq:tr} if $\Delta = \infty$.
	
	Let $\Delta < \infty$. If $\|d(0)\|_{L^1} \le \Delta$, then $d(0)$ minimizes \eqref{eq:tr}. We restrict to 
	$\|d(0)\|_{L^1} > \Delta$. Because  of the greedy construction, $d(\bar{x})$ is optimal
	if $\|d(\bar{x})\|_{L^1} = \Delta$ for some $\bar{x} < 0$. We consider the case where there is no such $\bar{x} \le 0$.
	We consider $\bar{x} \coloneqq \sup \{x \,|\, \|d(x)\|_{L^1} \le \Delta\}$.
	If $d$  is continuous at $\bar{x}$, then $\|d(\bar{x})\|_{L^1} = \Delta$, and $d(\bar{x})$ minimizes \eqref{eq:tr}.
	We distinguish two situations.
	
	\textbf{Situation 1.} If $d$ is only left continuous at $\bar{x}$, then we have
	$\|d(\bar{x})\|_{L^1} \le \Delta < \lim_{y\downarrow \bar{x}} \|d(y)\|_{L^1}$. Thus
	there exists a set $A \in \B$ satisfying $A \cap D(\bar{x}) = \emptyset$, $A \subset D(y)$ for all $y > \bar{x}$, and
	$\lambda(A) = \lim_{y\downarrow \bar{x}} \|d(y)\|_{L^1} - \|d(\bar{x})\|_{L^1}$.
	Such a set also exists if $\|d(\bar{x})\|_{L^1} \le \Delta$ and $d$ is neither left nor right continuous at $\bar{x}$.
	
	\textbf{Situation 2.} If $d$ is only right continuous at $\bar{x}$, then we have
	$\lim_{y\uparrow \bar{x}} \|d(y)\|_{L^1} \le \Delta < \|d(\bar{x})\|_{L^1}$.
	Thus there exists a set $A \in \B$ satisfying $A \subset D(\bar{x})$, $A \cap D(y) = \emptyset$ for all $y < \bar{x}$, and
	$\lambda(A) = \|d(\bar{x})\|_{L^1} - \lim_{y\downarrow \bar{x}} \|d(y)\|_{L^1}$.
	Such a set also exists if $ \Delta < \|d(\bar{x})\|_{L^1} $ and $d$ is neither left nor right continuous at $\bar{x}$.
	
	Because of the monotony of $\|d(\cdot)\|_{L^1}$ and the fact that the limits $\lim_{y\uparrow \bar{x}} \|d(y)\|_{L^1}$
	and $\lim_{y\downarrow \bar{x}} \|d(y)\|_{L^1}$ always exist by virtue of continuity from below and above of the Lebesgue
	measure, this distinction is exhaustive.
	Because of the mean value property of the Lebesgue measure \cite[Cor.\ 1.12.10]{bogachev2007measure},
	there exists $\B \ni B \subset A$ with $\lambda(B) = \Delta - \|d(\bar{x})\|_{L^1}$ (Situation 1)
	or $\lambda(B) = \|d(\bar{x})\|_{L^1} - \Delta$ (Situation 2). In Situation 1 we set
	$\tilde{D}\coloneqq D(\bar{x}) \cup B$, and in Situation 2 we set $\tilde{D} \coloneqq D(\bar{x}) \setminus B$.
	In both situations, $\tilde d \coloneqq \chi_{U\Delta \tilde{D}} -\chi_U$ minimizes \eqref{eq:tr}.
	
	The greedy construction of $D(\bar{x})$ and thus $\tilde{D}$ with respect to $g$  imply
	\begin{align*}
	\frac{1}{\Delta} \int_{\tilde{D}} g \dd s
	&\le \frac{1}{|\lambda(g^{-1}((-\infty,0)))|} \int_{g^{-1}((-\infty,0))} g \dd s\\
	&= -\frac{1}{|\lambda(g^{-1}((-\infty,0)))|} C(\chi_U) 
	\le -\frac{1}{|\lambda(\Omega)|} C(\chi_U),
	\end{align*}
	where we have used the definition of $C$ for the equality
	and that the integrand is negative for the second inequality.
	Then the identity $(\nabla J(\chi_U), \tilde{d})_{L^2} = \int_{\tilde{D}} g \dd s$ yields the claimed inequality.
\end{proof}

\section{Relationship between Set-Based and Characteristic Function Points of View}\label{sec:relationship}

In this appendix, we discuss the
relationship of \cref{ass:tr_convergence}
to \cite{hahn2022binary} and the corresponding
Taylor expansions.

\subsection{Relation of the setting of
	\cref{ass:tr_convergence} to \cite{hahn2022binary}}
Our setting and the assumptions cannot be compared
or embedded directly into the setting of
\cite{hahn2022binary} because \cite{hahn2022binary}
analyzes operations on sets in
atomless measure spaces; see, for example,
\cite[Definition\,1.12.7]{bogachev2007measure},
while we restrict \cref{alg:binary_tr_steepest_descent}
to functionals that operate on functions.
We note, however, that our arguments do not hinge on
the particular choice of the Lebesgue--Borel measure
for $L^1$ and $L^2$ and it is still possible to find correspondences of
the parts of \cref{ass:tr_convergence} in 
\cite{hahn2022binary}, which we do below.

First, we note that the $L^1$-regularization that is used in
the experiments in \cite{hahn2022binary} does satisfy
our assumptions because it is linear when
restricted to the feasible set of \eqref{eq:r}.

\Cref{ass:tr_convergence}, \ref{itm:Jboundedbelow}
is assumed in \cite[Theorem\,3]{hahn2022binary},
which essentially shows
$\liminf_{n\to\infty} C(\chi_{U^n}) = 0$.
We note that the assumption is not explicitly required if
\cref{ass:tr_convergence}, \ref{itm:Jprime_cc} holds
as well because the latter implies weak continuity of $J$
by means of \cref{prp:J_weakcont},
the feasible set of \eqref{eq:r} is weakly compact,
and continuous functions assume their minimum on
compact sets.

\Cref{ass:tr_convergence}, \ref{itm:Jdifferentiable} 
implies Assumption 1.3 in \cite{hahn2022binary}.
The local first-order Taylor expansion
in \cite[Theorem\,1]{hahn2022binary} for objective functions defined on measurable sets 
follows for the natural construction
\[ J_s(A) \coloneqq J(\chi_A), \text{ and } J_s'(A)D \coloneqq
\langle J'(\chi_A), \chi_{D\setminus A} - \chi_{D\cap A}
\rangle_{L^\infty(\Omega),L^1(\Omega)} \]
for $A$, $D \in \B$. We give a short proof in \cref{prp:taylor_relationship} below.

\Cref{ass:tr_convergence}, \ref{itm:Jprime_Lipschitz}
implies the assumption (5) in Lemma 3 and the (10) in
Theorem 3 in \cite{hahn2022binary}. They serve to
obtain sufficient decrease in the algorithm, which
is exactly the case, where it is needed the proof
of \cref{thm:asymptotics_part_1}
in this work. It is also similar
to Assumption 4.1 in \cite{leyffer2022sequential},
where it serves the same purpose.
It is no coincidence that assumptions of this type
are required for the analysis of descent algorithms that 
manipulate binary control functions in
$L^p$-norms for the following reason. All binary control functions $v$ satisfy
$\|v\|_{L^1} = \|v\|_{L^2}^2$. Therefore, bounding the error term of the Taylor
expansion by the squared $L^2$-norm is not sufficient to obtain a sufficient decrease 
condition because the linear predicted reduction is bounded from below
only by a fraction of the maximal $L^1$-norm of the trust-region step.
Thus we cannot prove that the linear predicted reduction dominates the remainder terms for
small trust-region radii without this further assumption.
Because the trust-region subproblem does not allow
fractional-valued control functions, a greedy strategy can always be used to approximate
the infimal value of the trust-region subproblem regardless of the $L^p$-norm
($p \in [0,\infty)$) that is used for the trust-region radius. Consequently, this
assumption cannot be avoided by choosing a different $L^p$-norm for the trust-region radius.

\Cref{ass:tr_convergence}, \ref{itm:Jprime_cc} is a compactness assumption on
the derivative of the objective function, which allows
us to infer the stationarity of weak-$^*$ cluster points
for \eqref{eq:r}.
It is not assumed in \cite{hahn2022binary}, which
does neither analyze the relationship to the continuous
relaxation nor show such a result.
It implies Assumption 1.4 for CIA in 
\cite{kirches2020compactness} by means of
\cref{prp:J_weakcont}.
The reason for this difference is that we need to pass to the
limit in the derivative of the objective functional in the norm when certifying stationarity.

\subsection{Taylor expansion for sets and characteristic functions}

Let $J_s$, $J'_s$ be given as above.
We say that $J_s$ is Fr\'{e}chet differentiable if
\[ J_s(A \symdif D)
= J_s(A) + J_s'(A)D + o(\lambda(D)), \]
which is the assertion of \cite[Theorem\,1]{hahn2022binary}.
Due to the assumed differentiability $J : L^1(\Omega) \to \R$
in \cref{ass:tr_convergence}, \ref{itm:Jdifferentiable},
that is with respect to the $L^1$-norm on the domain, we obtain that
$J_s$ is Fr\'{e}chet differentiable below.

\begin{proposition}\label{prp:taylor_relationship}
	$J_s$ is Fr\'{e}chet differentiable.
\end{proposition}
\begin{proof}
	Let $A$, $D \in B$. We use the defined identity $J_s(A\symdif D) = J(\chi_{A \symdif D})$,
	Taylor's theorem for $J$, and the identities
	$A \symdif D = (A \setminus D) \cup (D \setminus A)$,
	and $A = (A\setminus D) \cup (A \cap D)$---where both unions are disjoint---to deduce
	\begin{align*}
	J_s(A \symdif D) &= 
	J(\chi_A) + \langle J'(\chi_A), \chi_{A \symdif D} - \chi_A \rangle_{L^\infty,L^1}
	+ o(\|\chi_{A \symdif B} - \chi_A\|_{L^1})\\
	&= J(\chi_A) + \langle J'(\chi_A), \chi_{D\setminus A} - \chi_{A \cap D}\rangle_{L^\infty,L^1}
	+ o(\|\chi_{D\setminus A} - \chi_{A \cap D}\|_{L^1}).
	\end{align*}
	We observe that
	$\|\chi_{D\setminus A} - \chi_{A \cap D}\|_{L^1}
	= \|\chi_D\|_{L^1}$,
	and $\|\chi_D\|_{L^1} = \lambda(D)$. This implies
	\begin{gather*}
	J_s(A \symdif D) =
	J(\chi_A) + \langle J'(\chi_A), \chi_{D\setminus A} - \chi_{A \cap D}\rangle_{L^\infty,L^1}
	+ o(\lambda(D)).
	\end{gather*}
	Inserting the definitions of $J_s(A)$ and $J_s'(A)D$ yields the claim.
\end{proof}

\bibliographystyle{jnsao}
\bibliography{main}

\begin{thebibliography}{10}

\bibitem{alber1998projected}
Y.\,I{.\nobreak\kern 0.33333em}Alber, A.\,N{.\nobreak\kern 0.33333em}Iusem, and
  M.\,V{.\nobreak\kern 0.33333em}Solodov, On the projected subgradient method
  for nonsmooth convex optimization in a {H}ilbert space, \emph{Mathematical
  Programming} 81 (1998),  23--35,
  \href{https://dx.doi.org/10.1007/BF01584842}{\nolinkurl{doi:10.1007/bf01584842}}.

\bibitem{bendsoe2004extensions}
M.\,P{.\nobreak\kern 0.33333em}Bends{\o}e and O{.\nobreak\kern
  0.33333em}Sigmund, Extensions and applications, in \emph{Topology
  Optimization}, Springer, 2004,  71--158,
  \href{https://dx.doi.org/10.1007/978-3-662-05086-6}{\nolinkurl{doi:10.1007/978-3-662-05086-6}}.

\bibitem{bestehorn2022combinatorial}
F{.\nobreak\kern 0.33333em}Bestehorn, \emph{Combinatorial Algorithms and
  Complexity of Rounding Problems Arising in Mixed-Integer Optimal Control},
  PhD thesis, Technical University of Braunschweig, 2022,
  \href{https://dx.doi.org/10.24355/dbbs.084-202203101114-0}{\nolinkurl{doi:10.24355/dbbs.084-202203101114-0}}.

\bibitem{bestehorn2019switching}
F{.\nobreak\kern 0.33333em}Bestehorn, C{.\nobreak\kern 0.33333em}Hansknecht,
  C{.\nobreak\kern 0.33333em}Kirches, and P{.\nobreak\kern 0.33333em}Manns, A
  switching cost aware rounding method for relaxations of mixed-integer optimal
  control problems, in \emph{2019 IEEE 58th Conference on Decision and Control
  (CDC)}, 2019,  7134--7139,
  \href{https://dx.doi.org/10.1109/CDC40024.2019.9030063}{\nolinkurl{doi:10.1109/cdc40024.2019.9030063}}.

\bibitem{bestehorn2021mixed}
F{.\nobreak\kern 0.33333em}Bestehorn, C{.\nobreak\kern 0.33333em}Hansknecht,
  C{.\nobreak\kern 0.33333em}Kirches, and P{.\nobreak\kern 0.33333em}Manns,
  Mixed-integer optimal control problems with switching costs: a shortest path
  approach, \emph{Mathematical Programming} 188 (2021),  621--652,
  \href{https://dx.doi.org/10.1007/s10107-020-01581-3}{\nolinkurl{doi:10.1007/s10107-020-01581-3}}.

\bibitem{bestehorn2021switching}
F{.\nobreak\kern 0.33333em}Bestehorn, C{.\nobreak\kern 0.33333em}Hansknecht,
  C{.\nobreak\kern 0.33333em}Kirches, and P{.\nobreak\kern 0.33333em}Manns,
  Switching cost aware rounding for relaxations of mixed-integer optimal
  control problems: the 2-D case, \emph{IEEE Control Systems Letters} 6 (2021),
   548--553,
  \href{https://dx.doi.org/10.1109/CDC40024.2019.9030063}{\nolinkurl{doi:10.1109/cdc40024.2019.9030063}}.

\bibitem{boccardo1989non}
L{.\nobreak\kern 0.33333em}Boccardo and T{.\nobreak\kern
  0.33333em}Gallou{\"e}t, Non-linear elliptic and parabolic equations involving
  measure data, \emph{Journal of Functional Analysis} 87 (1989),  149--169,
  \href{https://dx.doi.org/10.1016/0022-1236(89)90005-0}{\nolinkurl{doi:10.1016/0022-1236(89)90005-0}}.

\bibitem{bogachev2007measure}
V.\,I{.\nobreak\kern 0.33333em}Bogachev, \emph{Measure Theory}, volume~1,
  Springer, 2007,
  \href{https://dx.doi.org/10.1007/978-3-540-34514-5}{\nolinkurl{doi:10.1007/978-3-540-34514-5}}.

\bibitem{dunn1980convergence}
J.\,C{.\nobreak\kern 0.33333em}Dunn, Convergence rates for conditional gradient
  sequences generated by implicit step length rules, \emph{SIAM Journal on
  Control and Optimization} 18 (1980),  473--487,
  \href{https://dx.doi.org/10.1137/0318035}{\nolinkurl{doi:10.1137/0318035}}.

\bibitem{hahn2022binary}
M{.\nobreak\kern 0.33333em}Hahn, S{.\nobreak\kern 0.33333em}Leyffer, and
  S{.\nobreak\kern 0.33333em}Sager, Binary optimal control by trust-region
  steepest descent, \emph{Mathematical Programming} 197 (2023),  147--190,
  \href{https://dx.doi.org/10.1007/s10107-021-01733-z}{\nolinkurl{doi:10.1007/s10107-021-01733-z}}.

\bibitem{haller2009holder}
R{.\nobreak\kern 0.33333em}Haller-Dintelmann, C{.\nobreak\kern 0.33333em}Meyer,
  J{.\nobreak\kern 0.33333em}Rehberg, and A{.\nobreak\kern 0.33333em}Schiela,
  H{\"o}lder continuity and optimal control for nonsmooth elliptic problems,
  \emph{Applied Mathematics and Optimization} 60 (2009),  397--428,
  \href{https://dx.doi.org/10.1007/s00245-009-9077-x}{\nolinkurl{doi:10.1007/s00245-009-9077-x}}.

\bibitem{hante2020mixed}
F.\,M{.\nobreak\kern 0.33333em}Hante, Mixed-integer optimal control for {PDE}s:
  relaxation via differential inclusions and applications to gas network
  optimization, in \emph{Mathematical Modelling, Optimization, Analytic and
  Numerical Solutions}, Springer, 2020,  157--171,
  \href{https://dx.doi.org/10.1007/978-981-15-0928-5_7}{\nolinkurl{doi:10.1007/978-981-15-0928-5_7}}.

\bibitem{hante2017challenges}
F.\,M{.\nobreak\kern 0.33333em}Hante, G{.\nobreak\kern 0.33333em}Leugering,
  A{.\nobreak\kern 0.33333em}Martin, L{.\nobreak\kern 0.33333em}Schewe, and
  M{.\nobreak\kern 0.33333em}Schmidt, Challenges in optimal control problems
  for gas and fluid flow in networks of pipes and canals: from modeling to
  industrial applications, in \emph{Industrial Mathematics and Complex
  Systems}, Springer, 2017,  77--122,
  \href{https://dx.doi.org/10.1007/978-981-10-3758-0_5}{\nolinkurl{doi:10.1007/978-981-10-3758-0_5}}.

\bibitem{hante2013relaxation}
F.\,M{.\nobreak\kern 0.33333em}Hante and S{.\nobreak\kern 0.33333em}Sager,
  Relaxation methods for mixed-integer optimal control of partial differential
  equations, \emph{Computational Optimization and Applications} 55 (2013),
  197--225,
  \href{https://dx.doi.org/10.1007/s10589-012-9518-3}{\nolinkurl{doi:10.1007/s10589-012-9518-3}}.

\bibitem{haslinger2015topology}
J{.\nobreak\kern 0.33333em}Haslinger and R.\,A{.\nobreak\kern
  0.33333em}M{\"a}kinen, On a topology optimization problem governed by
  two-dimensional {Helmholtz} equation, \emph{Computational Optimization and
  Applications} 62 (2015),  517--544,
  \href{https://dx.doi.org/10.1007/s10589-015-9746-4}{\nolinkurl{doi:10.1007/s10589-015-9746-4}}.

\bibitem{hinze2008optimization}
M{.\nobreak\kern 0.33333em}Hinze, R{.\nobreak\kern 0.33333em}Pinnau,
  M{.\nobreak\kern 0.33333em}Ulbrich, and S{.\nobreak\kern 0.33333em}Ulbrich,
  \emph{Optimization with {PDE} Constraints}, volume~23, Springer Science \&
  Business Media, 2008,
  \href{https://dx.doi.org/10.1007/978-1-4020-8839-1}{\nolinkurl{doi:10.1007/978-1-4020-8839-1}}.

\bibitem{jung2014relaxations}
M{.\nobreak\kern 0.33333em}Jung, \emph{Relaxations and Approximations for
  Mixed-integer Optimal Control}, PhD thesis, Heidelberg University, 2014,
  \href{https://dx.doi.org/10.11588/heidok.00016036}{\nolinkurl{doi:10.11588/heidok.00016036}}.

\bibitem{jung2015lagrangian}
M.\,N{.\nobreak\kern 0.33333em}Jung, G{.\nobreak\kern 0.33333em}Reinelt, and
  S{.\nobreak\kern 0.33333em}Sager, The Lagrangian relaxation for the
  combinatorial integral approximation problem, \emph{Optimization Methods and
  Software} 30 (2015),  54--80,
  \href{https://dx.doi.org/10.1080/10556788.2014.890196}{\nolinkurl{doi:10.1080/10556788.2014.890196}}.

\bibitem{kirches2020approximation}
C{.\nobreak\kern 0.33333em}Kirches, F{.\nobreak\kern 0.33333em}Lenders, and
  P{.\nobreak\kern 0.33333em}Manns, Approximation properties and tight bounds
  for constrained mixed-integer optimal control, \emph{SIAM Journal on Control
  and Optimization} 58 (2020),  1371--1402,
  \href{https://dx.doi.org/10.1137/18M1182917}{\nolinkurl{doi:10.1137/18m1182917}}.

\bibitem{kirches2020compactness}
C{.\nobreak\kern 0.33333em}Kirches, P{.\nobreak\kern 0.33333em}Manns, and
  S{.\nobreak\kern 0.33333em}Ulbrich, Compactness and convergence rates in the
  combinatorial integral approximation decomposition, \emph{Mathematical
  Programming} 188 (2021),  569--598,
  \href{https://dx.doi.org/10.1007/s10107-020-01598-8}{\nolinkurl{doi:10.1007/s10107-020-01598-8}}.

\bibitem{larsson1994class}
T{.\nobreak\kern 0.33333em}Larsson and M{.\nobreak\kern 0.33333em}Patriksson, A
  class of gap functions for variational inequalities, \emph{Mathematical
  Programming} 64 (1994),  53--79,
  \href{https://dx.doi.org/10.1007/BF01582565}{\nolinkurl{doi:10.1007/bf01582565}}.

\bibitem{leyffer2022sequential}
S{.\nobreak\kern 0.33333em}Leyffer and P{.\nobreak\kern 0.33333em}Manns,
  Sequential linear integer programming for integer optimal control with total
  variation regularization, \emph{ESAIM: Control, Optimisation and Calculus of
  Variations} 28 (2022), ~66,
  \href{https://dx.doi.org/10.1051/cocv/2022059}{\nolinkurl{doi:10.1051/cocv/2022059}}.

\bibitem{leyffer2021convergence}
S{.\nobreak\kern 0.33333em}Leyffer, P{.\nobreak\kern 0.33333em}Manns, and
  M{.\nobreak\kern 0.33333em}Winckler, Convergence of sum-up rounding schemes
  for cloaking problems governed by the {Helmholtz} equation,
  \emph{Computational Optimization and Applications} 79 (2021),  193--221,
  \href{https://dx.doi.org/10.1007/s10589-020-00262-3}{\nolinkurl{doi:10.1007/s10589-020-00262-3}}.

\bibitem{lindenstrauss1966short}
J{.\nobreak\kern 0.33333em}Lindenstrauss, A short proof of {L}iapounoff's
  convexity theorem, \emph{Journal of Mathematics and Mechanics} 15 (1966),
  971--972.

\bibitem{lyapunov1940completely}
A.\,A{.\nobreak\kern 0.33333em}Lyapunov, On completely additive vector
  functions, \emph{Izv. Akad. Nauk SSSR} 4 (1940),  465--478.

\bibitem{manns2021relaxed}
P{.\nobreak\kern 0.33333em}Manns, Relaxed multibang regularization for the
  combinatorial integral approximation, \emph{SIAM Journal on Control and
  Optimization} 59 (2021),  2645--2668,
  \href{https://dx.doi.org/10.1137/20M1377187}{\nolinkurl{doi:10.1137/20m1377187}}.

\bibitem{manns2020improved}
P{.\nobreak\kern 0.33333em}Manns and C{.\nobreak\kern 0.33333em}Kirches,
  Improved regularity assumptions for partial outer convexification of
  mixed-integer {PDE}-constrained optimization problems, \emph{ESAIM: Control,
  Optimisation and Calculus of Variations} 26 (2020), ~32,
  \href{https://dx.doi.org/10.1051/cocv/2019016}{\nolinkurl{doi:10.1051/cocv/2019016}}.

\bibitem{manns2020multidimensional}
P{.\nobreak\kern 0.33333em}Manns and C{.\nobreak\kern 0.33333em}Kirches,
  Multidimensional sum-up rounding for elliptic control systems, \emph{SIAM
  Journal on Numerical Analysis} 58 (2020),  3427--3447,
  \href{https://dx.doi.org/10.1137/19M1260682}{\nolinkurl{doi:10.1137/19m1260682}}.

\bibitem{manns2021approximation}
P{.\nobreak\kern 0.33333em}Manns, C{.\nobreak\kern 0.33333em}Kirches, and
  F{.\nobreak\kern 0.33333em}Lenders, Approximation properties of sum-up
  rounding in the presence of vanishing constraints, \emph{Mathematics of
  Computation} 90 (2021),  1263--1296.

\bibitem{sager2005numerical}
S{.\nobreak\kern 0.33333em}Sager, \emph{Numerical Methods for Mixed-integer
  Optimal Control Problems}, Der Andere Verlag L{\"u}beck, 2005.

\bibitem{sager2012integer}
S{.\nobreak\kern 0.33333em}Sager, H.\,G{.\nobreak\kern 0.33333em}Bock, and
  M{.\nobreak\kern 0.33333em}Diehl, The integer approximation error in
  mixed-integer optimal control, \emph{Mathematical Programming} 133 (2012),
  1--23,
  \href{https://dx.doi.org/10.1007/s10107-010-0405-3}{\nolinkurl{doi:10.1007/s10107-010-0405-3}}.

\bibitem{sager2011combinatorial}
S{.\nobreak\kern 0.33333em}Sager, M{.\nobreak\kern 0.33333em}Jung, and
  C{.\nobreak\kern 0.33333em}Kirches, Combinatorial integral approximation,
  \emph{Mathematical Methods of Operations Research} 73 (2011),  363--380,
  \href{https://dx.doi.org/10.1007/s00186-011-0355-4}{\nolinkurl{doi:10.1007/s00186-011-0355-4}}.

\bibitem{sharma2020inversion}
M{.\nobreak\kern 0.33333em}Sharma, M{.\nobreak\kern 0.33333em}Hahn,
  S{.\nobreak\kern 0.33333em}Leyffer, L{.\nobreak\kern 0.33333em}Ruthotto, and
  B{.\nobreak\kern 0.33333em}van Bloemen~Waanders, Inversion of
  convection--diffusion equation with discrete sources, \emph{Optimization and
  Engineering} 22 (2021),  1419--1457,
  \href{https://dx.doi.org/10.1007/s11081-020-09536-5}{\nolinkurl{doi:10.1007/s11081-020-09536-5}}.

\bibitem{tartar1979compensated}
L{.\nobreak\kern 0.33333em}Tartar, Compensated compactness and applications to
  partial differential equations, in \emph{Nonlinear Analysis and Mechanics:
  {H}eriot--{W}att Symposium}, volume~4, 1979,  136--212.

\bibitem{toint1997non}
P.\,L{.\nobreak\kern 0.33333em}Toint, Non-monotone trust-region algorithms for
  nonlinear optimization subject to convex constraints, \emph{Mathematical
  Programming} 77 (1997),  69--94,
  \href{https://dx.doi.org/10.1007/bf02614518}{\nolinkurl{doi:10.1007/bf02614518}}.

\bibitem{vogt2020mixed}
R.\,H{.\nobreak\kern 0.33333em}Vogt, S{.\nobreak\kern 0.33333em}Leyffer, and
  T{.\nobreak\kern 0.33333em}Munson, A mixed-integer {PDE}-constrained
  optimization formulation for electromagnetic cloaking, \emph{SIAM Journal on
  Scientific Computing} 44 (2022),  B29--B50,
  \href{https://dx.doi.org/10.1137/20m1315993}{\nolinkurl{doi:10.1137/20m1315993}}.

\bibitem{yu2021multidimensional}
J{.\nobreak\kern 0.33333em}Yu and M{.\nobreak\kern 0.33333em}Anitescu,
  Multidimensional sum-up rounding for integer programming in optimal
  experimental design, \emph{Mathematical Programming} 185 (2021),  37--76,
  \href{https://dx.doi.org/10.1007/s10107-019-01421-z}{\nolinkurl{doi:10.1007/s10107-019-01421-z}}.

\bibitem{zeile2021mixed}
C{.\nobreak\kern 0.33333em}Zeile, N{.\nobreak\kern 0.33333em}Robuschi, and
  S{.\nobreak\kern 0.33333em}Sager, Mixed-integer optimal control under minimum
  dwell time constraints, \emph{Mathematical Programming} 188 (2021),
  653--694,
  \href{https://dx.doi.org/10.1007/s10107-020-01533-x}{\nolinkurl{doi:10.1007/s10107-020-01533-x}}.

\end{thebibliography}

\vfill

\framebox{\parbox{.92\linewidth}{The submitted manuscript has been created by
		UChicago Argonne, LLC, Operator of Argonne National Laboratory (``Argonne'').
		Argonne, a U.S.\ Department of Energy Office of Science laboratory, is operated
		under Contract No.\ DE-AC02-06CH11357.  The U.S.\ Government retains for itself,
		and others acting on its behalf, a paid-up nonexclusive, irrevocable worldwide
		license in said article to reproduce, prepare derivative works, distribute
		copies to the public, and perform publicly and display publicly, by or on
		behalf of the Government.  The Department of Energy will provide public access
		to these results of federally sponsored research in accordance with the DOE
		Public Access Plan \url{http://energy.gov/downloads/doe-public-access-plan}.}}

\end{document}